\documentclass[12pt]{amsart}
\usepackage{amsmath,amsfonts,amssymb,amscd}
\title{An energy-theoretic approach to the\\
Hitchin-Kobayashi correspondence  \\
for manifolds, II}
\author{Toshiki Mabuchi${}^*$}
\address{Department of Mathematics, Graduate School of Science, Osaka University, 
Toyonaka, Osaka, 560-0043 Japan}
\date{}
\begin{document}
\maketitle

\thanks{\em \qquad Dedicated to Professor Eugenio Calabi on his eightieth birthday}

\abstract  Recently, Donaldson proved  asymptotic stability for a polarized algebraic manifold $M$
with  polarization class admitting a K\"ahler metric of constant scalar curvature, essentially when the
linear algebraic part $H$ of $\operatorname{Aut}^0(M)$ is semisimple.
The purpose of this paper is to give a generalization of Donaldson's result to the case 
where the polarization class
admits an extremal K\"ahler metric, even when 
$H$ is not semisimple.
\endabstract

\footnotetext{\; \; ${}^*$Partially supported 
by JSPS Grant-in-Aid for Scientific Research B No.~13440023.}

\section{Introduction}
For  a connected polarized algebraic manifold $(M, L)$ with an extremal K\"ahler metric
in the polarization class $c_1(L)_{\Bbb R}$, we consider the Kodaira
embedding 
$$
\Phi_m = \Phi_{|L_{}^{\otimes m}|} : M\; \hookrightarrow \;\Bbb P^*(V_m), \qquad m \gg 1,
$$
where $\Bbb P^*(V_m)$ denotes the set of all hyperplanes in 
$V_m := H^0(M, \mathcal{O}(L^m))$ through the origin.
For the identity component $\operatorname{Aut}^0(M)$ of the 
group of holomorphic automorphisms of $M$,
let $H$ denote its maximal connected
linear algebraic subgroup.
Replacing the ample holomorphic  line bundle $L$ 
by some positive integral multiple of $L$ if necessary, 
we may assume that the 
natural $H$-equivariant maps
$$
\operatorname{pr}_m  : \; \otimes^m V_1 \to\, V_m, \qquad m = 1,2, \dots, 
$$
are surjective (cf. \cite{M}, Theorem 3), and 
may further fix an $H$-linearization of $L$, i.e., a lift to $L$ of the $H$-action on $M$ 
such that $H$ acts on $L$ as bundle isomorphisms covering the $H$-action on $M$.
In this paper, applying a method in \cite{M5}, we shall generalize a result 
in Donaldson \cite{D}  about stability to extremal K\"ahler cases:

\medskip\noindent
{\bf Main Theorem.}  \,
{\em For a polarized algebraic manifold $(M, L)$ as above with an extremal K\"ahler
metric in the polarization class, there exists
an algebraic torus $T$ in $H$
such that the image $\Phi_m (M)$ in $\Bbb P^*(V_m)$ is stable relative to $T$
{\rm (cf. Section 2 and \cite{M3})} for $m \gg 1$.}

\medskip
In particular in \cite{M4}, by an argument as in \cite{D}, 
an extremal K\"ahler metric in a fixed integral K\"ahler class on a projective algebraic manifold
$M$  will be shown to be unique up to the action of the group $H$.

\medskip
We now consider the set $\Delta_{L}$ of all algebraic tori $T$ in $H$
for which the statement of Main Theorem is valid. 
Note that, if $T'$ and  $T''$ are algebraic tori in $H$
with $T' \subset T''$ and  $T' \in \Delta_{L}$, then $T''$ also belongs to
$\Delta_{L}$ (cf. \cite{M3}, Theorem 3.2).  Fix once for all an extremal K\"ahler metric $\omega_0$ in the
polarization class in Main Theorem.
By a result of Calabi \cite{Ca},
the identity component $K$ of the group of isometries of $(M, \omega_0 )$
is a maximal compact connected subgroup of $H$. 
For the identity component $Z$ of the center of $K$, 
we consider the complexification $Z^{\Bbb C}$ of $Z$ in  $H$. 
In order to prove Main Theorem, it suffices to show

\medskip\noindent
{\bf Theorem I.}  \;
{\em The torus $Z^{\Bbb C}$ belongs to $\Delta_{L}$.}

\medskip
For the set $\Delta$ of all algebraic tori $T$ in $Z^{\Bbb C}$,
we consider its subset $\Delta_0$
consisting of all irredundant elements in $\Delta$
(cf. Definition 2.4), where $Z^{\Bbb C}$ always belongs to $\Delta_0$.
Furthermore, related to the obstruction as in \cite{M2}, we shall show that
$$
Z^{\Bbb C}\in \Delta_{1},
\leqno{(1.1)}
$$
where $\Delta_1$ is the subset of $\Delta$ as in Section 2 below, and is obtained
from a weighted version of Donaldson's asymptotic
expansion approximating balanced metrics.
Then the proof of Theorem I is  reduced to showing the 
following:

\medskip\noindent
{\bf Theorem II.}  \;\; 
$\Delta_{L} \cap \Delta_0 \; = \; \Delta_0\cap \Delta_1$.

\medskip
If the scalar curvature of the extremal K\"ahler metric $\omega_0$ above is
constant, and if the obstruction as in \cite{M2} vanishes, 
then we have both
$\Delta_1 =\Delta$
 and
$\{1\} \in \Delta_0$.
 Hence in this case, Theorem II shows that
$\{1\}$ sits in $\Delta_{L}$.
This then proves the main theorem in \cite{M5}. 

\medskip
One may ask why relative stability in place of ordinary stability has to be considered in our study.
The reason why we choose relative stability is because, in general, the obstruction in \cite{M2}
to asymptotic semistability does not vanish (cf. \cite{MN}). Thus, as to the group action on $V_m$ related
to stability, we must replace the full special linear group
$\operatorname{SL}(V_m)$ of $V_m$ by its subgroup $G_m(T)$ (see (2.3)), where 
the algebraic torus $T$ in $Z^{\Bbb C}$ is chosen in such a way that the obstruction vanishes for the group
$G_m (T)$, i.e., $G_m (T)$ fixes $\hat{M}_m$
(cf. Section 2). Note also that $G_m (T)$ is a direct product of special linear
groups.
To see why we choose such a group $G_m (T)$ in place of $\operatorname{SL}(V_m)$, we compare our stability 
with that of holomorphic vector bundles. Recall
that a holomorphic vector bundle splitting into a
direct sum of stable vector bundles often appears in the boundary of a compactified moduli space of stable
vector bundles. Similarly for our stability of manifolds, a splitting
phenomenon occurs for
$V_m$ in (2.2). Roughly speaking, we consider the moduli space of all $M$'s with fixed
decomposition  data (2.2), where same
type of construction of moduli spaces occurs typically for the Hodge decomposition in the variation of Hodge
structures.

\medskip
We now explain the difficulty which we encounter in applying the method of \cite{M5}.
Such a difficulty comes up when we use the estimate of Phong and Sturm \cite{PS}.
By applying a stability criterion in \cite{M5} of Hilbert-Mumford's  type,
we write the vector space $\frak p_m$ as an orthogonal direct sum
 $$
\frak p_m \;  =\;  \frak p_m' \oplus\frak p_m'', \;\; \quad\text{(cf. Section 3),}
$$
and then check the stability of $\hat{M}_m$ along the orbits of the one-parameter subgroups in $G_m
(T)$ generated by elements of
$\frak p_m''$. Though $\frak p$ and $\frak p_m''$ are transversal by 
the equality $\frak p'_m = \frak p_m \cap \frak p$, 
we further need the orthogonality of $\frak p$ and $\frak p''_m$ in order to apply directly 
the estimate in \cite{PS}. 
Since such an orthogonality does not generally hold, we are in trouble, 
but still the situation is not so bad (see
(3.17), (3.18)), and this overcomes the difficulty.

\medskip
It is very likely that the set $\Delta_L$ has a natural minimal element closely related to the construction of
the torus $T_0$ in Section 2. To see this, let us consider the case where $M$
is an extremal K\"ahler toric Fano surface polarized by $L = K_M^{-1}$.  Then $M$ is possibly a 
complex projective plane blown up at $r$ points with $r \leq 3$. If $r = 0$ or $3$, then $M$
admits a K\"ahler-Einstein metric, and the explanation following Theorem II above shows that
$\Delta_L$ has the unique minimal element $\{1\} \, (=\, T_0)$. On the other hand, if $r =1$, then 
$T_0$ coincides with $Z^{\Bbb C}$, and is the one-dimensional torus generated by the extremal K\"ahler
vector field. Hence, in this case, $\Delta_L$ has the unique minimal element $T_0$.
Finally for $r = 2$, the involutive holkomorphic symmetry of $M$ switching the blown-up points allows us to
regard $T_0$ as the one-dimensional torus generated by the extremal K\"ahler vector field.
It then follows that $T_0$ again has to be a minimal element of $\Delta_L$.

\section{Notation, convention and preliminaries}

\smallskip
Throughout this paper, we fix once for all a pair $(M, L)$ of a connected projective algebraic manifold $M$
and an ample holomorphic line bundle $L$ over $M$ as in the
introduction. For $V_m$ in the introduction, we put $N_m := \dim_{\Bbb C} V_m -1$, where 
the positive integer $m$ is such that $L^m$ is
very ample. Let $n$ and $d$ be respectively the dimension of $M$ and the
degree of  the image $M_m := \Phi_ m
(M )$ in the projective space $\Bbb P^* (V_m)$. 
Fixing an $H$-linearization of $L$ as in the introduction, 
we consider the associated representation: $H \to \operatorname{PGL}(V_m)$.
Pulling it back by
the finite unramified cover: $\operatorname{SL}(V_m) \to \operatorname{PGL}(V_m)$,
we obtain an isogeny
$$
\iota : \tilde{H} \to H,
\leqno{(2.1)}
$$
where $\tilde{H}$ is an algebraic subgroup of $\operatorname{SL}(V_m)$.
On the other hand, for an algebraic torus $T$ in $H$,
the $H$-linearization of $L$ naturally induces a faithful representation
$$
H \to \operatorname{GL}(V_m),
$$
and this gives a $T$-action on $V_m$ for each $m$.
Then we have a finite subset $\Gamma_m = \{ \chi_1, \chi_2, \dots, \chi_{\nu_m}\}$
 of the free Abelian group $\operatorname{Hom}(T, \Bbb C^*)$ of all characters of $T$
 such that the vector space
$V_m = H^0(M, \mathcal{O}(L^m))$ is uniquely written as a direct sum
$$
V_m\; =\; \bigoplus_{k=1}^{\nu_m}\; V_T(\chi_k ),
\leqno{(2.2)}
$$
where for each $\chi \in \operatorname{Hom}(T, \Bbb C^*)$, we set
$V_T(\chi ) :=
\{ s  \in V_m \, ;\, t \cdot s = \chi (t ) \, s$ for all 
$t \in T \}$.
Define an algebraic subgroup $G_m = G_m (T)$ of $\operatorname{SL}(V_m)$ by
$$
G_m\;  :=\;  \prod_{k=1}^{\nu_m}\;\operatorname{SL}(V_T(\chi_k )),
\leqno{(2.3)}
$$
and the associated Lie subalgebra of $\operatorname{sl}(V_m)$ 
will be denoted by $\frak g_m$.
Here, $G_m$ and $\frak g_m$ possibly
depend on the choice of the algebraic torus $T$,
and if necessary, we denote these by $G_m (T)$ and $\frak g_m (T)$, respectively.
The $T$-action on $V_m$ is, more precisely, a right action, while  
 the $G_m$-action on $V_m$
is a left action. Since $T$ is Abelian, this $T$-action on $V_m$ 
can be regarded also as a left action.
Note that the  group $G_m$ acts diagonally on
$V_m$  in such a way that, for each $k$,
the $k$-th factor $\operatorname{SL}(V_T(\chi_k ))$ of $G_m$
acts  just on the $k$-th factor 
$V_T(\chi_k )$ of $V_m$. 
We now put
$$
W_m :=\{S^d (V_m )\}^{\otimes n+1},
$$
where $S^d(V_m)$ denotes the $d$-th symmetric tensor product of $V_m$.
To the image $M_m$ of $M$,
we can associate a nonzero element
$\hat{M}_m$ in $W^*_m$ such that 
the corresponding element $[\hat{M}_m ]$ in 
$\Bbb P^* (W_m)$ is the Chow point of the 
irreducible reduced algebraic cycle
$M_m$ on $\Bbb P^* (V_m)$. 
Note that the $G_m$-action on $V_m$
naturally induces a $G_m$-action on $W_m$ and also on
$W_m^*$. As in \cite{M3}, the subvariety $M_m$ of $\Bbb P^*(V_m)$ is 
said to be {\it stable relative to $T$\/} or {\it semistable relative to $T$\/}, according as the orbit 
$G_m \cdot \hat{M}_m$ is closed in $W^*_m$ or the closure of $G_m \cdot \hat{M}_m$ in $W^*_m$
does not contain the origin of $W^*_m$.

\medskip
Take a Hermitian metric $h_0$ for $L$ such that
$c_1(L; h_0 )$ is the extremal K\"ahler metric $\omega_0$ in the Main Theorem.
Let $E$ be the extremal K\"ahler vector field for $(M, \omega_0
)$, and let $\frak k$ be the Lie algebra of $K$. 
For  $\omega_0$ above, 
we further define $\Delta_{\operatorname{min}}$ as the set of all $T\in \Delta$
for which the statement of Theorem B in \cite{M3} is valid.
Then, as the procedure in Section 6 of \cite{M3} shows, 
there exists a unique minimal element, denoted by $T_0$,
of $\Delta_{\operatorname{min}}$ such that
$\Delta_{\operatorname{min}} \, =\, \{\, T\in \Delta \, ;\, T_0 \subset T\,\}$.
Then by the notation in (3.1) below, 
$T_0$ is the closure in $Z^{\Bbb C}$ of the complex Lie
subgroup generated by the vector fields 
$$
E, \; F_k, \quad k =1, 2, \dots\,\,, 
$$
which appear in the asymptotic expansion  
approximating weighted analogues (cf. \cite{M3}, 2.6) of balanced metrics.
For each $T\in \Delta_{\operatorname{min}}$, we put 
$\tilde{T}\;  := \; \iota^{-1}(T)$,
and let $G'_m (T)$ and $Z_m'(T)$ be the identity components of  $G_m (T)\cap \tilde{H}$
and $G_m (T) \cap Z^{\Bbb C}$, respectively.  Put
$$
\Delta_1  := \, \{\, T\in \Delta_{\operatorname{min}} \,;\, 
G'_m (T)\cdot \hat{M}_m = \hat{M}_m\,
\}.
$$
{\em Definition $2.4$.} 
For an algebraic torus $T$ in $\Delta$, we say that $T$ is 
{\it irredundant}, if $\dim_{\Bbb C} K^{\Bbb C} = \dim_{\Bbb C} G_m'(T) + \dim_{\Bbb C} T$, 
or equivalently if $\dim_{\Bbb C} Z^{\Bbb C} = 
\dim_{\Bbb C} Z_m'(T) + \dim_{\Bbb C} T$.
For instance, if $T = Z^{\Bbb C}$, then $Z'_m(T) = \{1\}$.
Hence, $Z^{\Bbb C}$ is irredundant.

\medskip
For the complexification $K^{\Bbb C}$ of $K$ in $H$,
we define $\tilde{K}^{\Bbb C}:= \iota^{-1}(K^{\Bbb C})$ and its subset $\tilde{Z}^{\Bbb C}:= \iota^{-1}(Z^{\Bbb C})$.
For the time being, let $T = Z^{\Bbb C}$, and we consider the associated  
set $\Gamma_m = \{\chi_1, \chi_2. \dots, \chi_{\nu_m}\}$ of characters.
Since $K^{\Bbb C}$ commutes with $Z^{\Bbb C}$, 
we have the inclusion
$$
\tilde{K}^{\Bbb C} \;\;\;  \subset \;\; \;  \operatorname{SL}(V_m )\; \bigcap\;\;  
\prod_{k=1}^{\nu_m} \;
\operatorname{GL}(V^{}_{Z^{\Bbb C}}(\chi_k )).
\leqno{(2.5)}
$$
Recall that the extremal K\"ahler vector field $E$ belongs to the
Lie algebra of
$T_0$. Hence, a theorem of Calabi \cite{Ca} shows that 
$G'_m (Z^{\Bbb C}) \, \subset\, G'_m (T_0)\, \subset\,  \tilde{K}^{\Bbb C}$. Hence,
$$
G'_m (Z^{\Bbb C})\cdot \tilde{Z}^{\Bbb C}\, \subset\, \tilde{K}^{\Bbb C}.
\leqno{(2.6)}
$$

\smallskip\noindent
{\it Proof of $(1.1)$.}
To show (1.1), we compare two groups $[\tilde{K}^{\Bbb C}, \tilde{K}^{\Bbb C}]$ and
$G'_m (Z^{\Bbb C})$. By (2.5), we obviously have
$[\tilde{K}^{\Bbb C}, \tilde{K}^{\Bbb C}] \, \subset\, G'_m (Z^{\Bbb C})$. 
On the other hand,
$$
\dim_{\Bbb C}\, [\tilde{K}^{\Bbb C}, \tilde{K}^{\Bbb C}] \;
= \; \dim_{\Bbb C}\, \tilde{K}^{\Bbb C} - \dim_{\Bbb C}\, \tilde{Z}^{\Bbb C} \; \geq \; 
\dim_{\Bbb C}\, G'_m (Z^{\Bbb C}),
$$
where the last inequality follows from (2.6) in view of the fact that 
the intersection of $G'_m (Z^{\Bbb C})$
and $\tilde{Z}^{\Bbb C}$ is a finite group.
Now, we see that $G'_m (Z^{\Bbb C})$ coincides with $[\tilde{K}^{\Bbb C}, \tilde{K}^{\Bbb C}]$.
Hence $G'_m (Z^{\Bbb C})\cdot \hat{M}_m = \hat{M}_m$.
Then by $T_0 \subset Z^{\Bbb C}$, we now obtain $Z^{\Bbb C} \in \Delta_1$, as required.
\qed

\medskip
Let $h$ be  a Hermitian metric for $L$
such that $\omega := c_1(L; h)$ is a $K$-invariant K\"ahler metric on $M$.
Define a Hermitian metric $ \rho_h$ on $V_m$ by
$$
\rho_h (s, s')\;\; := \int_M (s, s')^{}_{h^m} \, \omega^n, \qquad s, s'\in V_m,
\leqno{(2.7)}
$$
where $(s, s')^{}_{h^m}$ denotes the function on $M$ 
obtained as the the pointwise inner product of $s$, $s'$ by $h^m$.
Let $\mathcal{S} := \{s_0, s_1, \dots, s_{N_m}\}$ 
be an orthonormal basis for
$V_m$ satisfying
$$
\rho_h (s_i, s_j)\;  =\; \delta_{ij}.
$$
Let $T \in \Delta_0 \cap \Delta_1$.
Then we say that $\mathcal{S}$ is $T$-{\it admissible}, if each $V_T(\chi_k )$,
$k =1,2, \dots$,
admits a basis $\{s_{k,i}\,;\, i =1,2,\dots,n_k\}$ such that
 $$
s_{l (k,i )} = s_{k,i}, \quad\qquad i = 1,2,\dots, n_k;\; \; k =1,2, \dots,
\nu_m,
\leqno{(2.8)}
$$
where $n_k := \dim_{\Bbb C}\, V_T(\chi_k )$, and
$l (k, i) \, :=\, (i - 1)\, +\, \Sigma_{k'=1}^{k-1} n_{k'}$
for all $k$ and $i$ (cf. \cite{M3}).
Let $\frak t_c := \operatorname{Lie}(T_c)$ denote the Lie algebra of the maximal compact
subgroup
$T_c$ of
$T$. Put $ q := 1/m$ and
$\frak t_{\Bbb R} := \sqrt{-1}\, \frak t_c$. For  each $F\in \frak t_{\Bbb R}$, we define
$$
B_q(\omega, F) \; :=\; \frac{n!}{m^n}\sum_{k =1}^{\nu_m}\sum_{i =1}^{n_k} \; 
e^{-d\chi^{}_k(F)}\, |s_{k,i}|^2_{h^m},
\leqno{(2.9)}
$$
where $|s|^2_{h^m} := (s, s)_{h^m}$
for all $s \in V_m$, and
$d\chi_k : \frak t_{\Bbb R} \to \Bbb R$ denotes the restriction to $\frak t_{\Bbb R}$ 
of the differential at $t =1$ for the character 
$\chi_k \in \operatorname{Hom}(T, \Bbb C^*)$.

\smallskip
As a final remark in this section, we give an upper bound for degrees of the characters in
$\Gamma_m$. Let $T$ be an algebraic torus sitting in $Z^{\Bbb C}$.
By setting $r := \dim_{\Bbb C} T$, we identify
$T$ with the multiplicative group $(\Bbb C^* )^r := \{t = (t_1, t_2, \dots, t_r ) \,;\, t_j \in \Bbb C^*
\text{ for all $j$}\}$. Since each $\chi_k$ in (2.2) may depend on $m$, the character $\chi_k$ 
will be rewritten as $\chi_{m;k}$ until the end of this section.
Then for each $k \in \{1,2,\dots, \nu_m\}$, 
$$
\chi^{}_{m;k} (t) \; =\; \prod_{i=1}^r \, t_i^{\alpha (m,k,i)},
\qquad t = (t_1, t_2, \dots. t_r) \in T,
$$
for some integers $\alpha (m,k,i)$ independent of the choice of $t$.
Define a nonnegative integer $\alpha_m$ by 
$\alpha_m := \sup_{k=1}^{\nu_m} \Sigma_{i=1}^r |\alpha (m,k,i)|$. 
Then we have the following upper bound for $\alpha_m$:

\medskip\noindent
{\bf Lemma 2.10.} \; 
{\em For all positive integer m, the inequality $\alpha_m\, \leq\, m\, \alpha_1$ holds.}

\medskip\noindent
{\em Proof}: 
Put $S: = \operatorname{Ker}\,  \operatorname{pr}_m$. Since the subspace $S$ of $\otimes^m V_1$ is preserved 
by the $T$-action, we have a $T$-invariant subspace, denoted by $S^{\perp}$,
of $\otimes^m V_1$ such that the vector space $\otimes^m V_1$ is written as  a direct sum
$$
\otimes^m V_1 \; =\; S \oplus S^{\perp}.
$$
Then the restriction of $\operatorname{pr}_m$ to $S^{\perp}$ defines a $T$-equivariant isomorphism
$S^{\perp} \cong V_m$. On the other hand, the characters of $T$ appearing in the $T$-action on $\otimes^m V_1$
are
$$
\;\chi^{}_{\vec{k}} (t) \, := \, 
 t_1^{\Sigma_{j=1}^m \alpha (1, k_j, 1)} t_2^{\Sigma_{j=1}^m\alpha (1, k_j, 2)} \cdots 
t_r^{\Sigma_{j=1}^m\alpha (1, k_j, r)},
\; \vec{k} = (k_1, k_2, \dots, k_m)\in I^m,
$$
where $I^m$ is the Cartesian product of $m$-pieces of $I := \{1,2,\dots, \nu_m\}$.
Since $S^{\perp} \, (\cong V_m )$ is a subspace of $\otimes^m V_1$, 
we now obtain 
$$
\alpha_m \;\leq \; \max_{\;\vec{k} \in I^m} \Sigma_{i=1}^r |\Sigma_{j=1}^m \,\alpha (1, k_j, i)|
\; \leq \; \max_{\;\vec{k} \in I^m}\Sigma_{j=1}^m\Sigma_{i=1}^r | \alpha (1, k_j, i)|
\; \leq \; m\, \alpha_1,
$$
as required.
\qed
\medskip

\noindent{\em Remark $2.11.$}\; 
By the definition of $F( \ell )$ in (3.1) below, Lemma 2.10 above implies that
$|d \chi_{m;k}( F( \ell ))| \leq C \alpha_1 q $ for some positive real constant $C$ independent of the choice of $m$ and $k$.
Hence in $(3.1)$ below, for each fixed nonnegative integer $\ell$, there exists a 
positive constant $C'$ independent of $m$ and $k$ such that
$$
| e^{-d\chi_k (F (\ell ))} - 1 | \; \leq \;  C' q, \qquad k = 1,2,\dots, \nu_m.
$$ 
In particular, 
by the notation in $(3.2)$ below, the integral
$\int_M \| \tilde{s}_{k,i} \|^2_{h^m} \omega (\ell )^n\, (= e^{-d\chi_k (F (\ell ))}) $
converges to $1$, uniformly in $k$, as $m\to\infty$.

\section{Proof of Theorem II}

\smallskip
Fix an arbitrary element $T$ of $\Delta_0\cap \Delta_1$. Let $m \gg 1$. 
Then by \cite{M3}, Theorem B, 
there exist  
$F_k\in \frak t_{\Bbb R}$, real numbers $\alpha_k \in \Bbb R$, 
and smooth real-valued $K$-invariant functions $\varphi_k$, $k = 1$,$2$,\dots,  on $M$ such
that, for each $\ell\in \Bbb Z_{\geq 0}$,  we have
$$
B_q (\omega (\ell ), F (\ell )) 
\; =\; C_{q,\ell} + 0(q^{\ell +2}),  \qquad m \gg 1,
\leqno{(3.1)}
$$
where
$F (\ell ) := (\sqrt{-1}\, E/2)\,q^2\, +\, \Sigma_{j=1}^{\ell}\, q^{j+2}
F_j$, 
$h(\ell ):= h_0 \exp (-\Sigma_{k=1}^{\ell} q^j\varphi_j )$, $C_{q, \ell} := 1 + \Sigma_{j=0}^{\ell}\, \alpha_j
q^{j+1}$, and $\omega (\ell ) := c_1(L; h(\ell ))$. 
Let us now fix an arbitrary positive integer $\ell$.
To each $T$-admissible orthonormal basis $\mathcal{S} := \{s_0, s_1, \dots, s_{N_m}\}$ 
for $(V_m; \rho_{h(\ell )})$, we associate a basis $\tilde{\mathcal{S}} :=
 \{\tilde{s}_0, \tilde{s}_1, \dots, \tilde{s}_{N_m}\}$ for $V_m$ by
$$
\tilde{s}_{k,i} \; =\; e_{}^{-d\chi_k (F(\ell ))/2}{s}_{k,i},
\quad\qquad i = 1,2,\dots, n_k;\; \; k =1,2, \dots,
\nu_m,
\leqno{(3.2)}
$$
where we put $s_{l (k, i ) } = s_{k, i}$ and $\tilde{s}_{l (k, i ) } = \tilde{s}_{k,i}$ by
using the notation in (2.8).
We now consider the Kodaira embedding $\Phi_m : M \to \Bbb P^* (V_m )$ defined by
$$
\Phi_m (x)\,  :=\,  (\tilde{s}_0 (x): \tilde{s}_1(x): \dots : \tilde{s}_{N_m}(x)),
\qquad x \in M,
$$
where $\Bbb P^* (V_m )$ is identified with $\Bbb P^{N_m}(\Bbb C ) = \{\,(z_0:z_1:\dots :z_{N_m})\,\}$ 
by the basis $\tilde{\mathcal{S}}$. Put $M_m := \Phi_m (M)$.
Since $\Delta_L \cap \Delta_0$ is a subset of  $\Delta_1 \cap \Delta_0$
(cf. \cite{M2}, Section 3), 
the proof of Theorem II (and Main Theorem also) is reduced to showing 
the following assertion:

\medskip\noindent
{\em Assertion$:$
\qquad The orbit $G_m (T) \cdot \hat{M}_m$ is closed in $W_m^*$.}

\medskip
In the Hermitian vector space $(V_m; \rho_{h(\ell )})$,
the subspaces  $V_T (\chi_k )$, $k$ =1,2,\dots, $\nu_m$, are mutually orthogonal. 
Put
$$
K_m \; :=\; \prod_{k=1}^{\nu_m}\;\operatorname{SU}(V_T (\chi_k ); \rho_{h(\ell )}),
\qquad \frak k_m \: :=\; \bigoplus_{k=1}^{\nu_m}\;\frak s\frak u (V_T (\chi_k ); \rho_{h(\ell
)}).
$$
Since $T$ belongs to $\Delta_1$, the group $G'_m = G'_m (T)$ coincides with 
 the isotropy subgroup of $G_m$ at $\hat{M}_m\in
W_m^*$.  Consider the Lie algebra $\frak g'_m := \operatorname{Lie}(G'_m)$ of $G'_m$.
Put $\frak h := \operatorname{Lie} (H) = 
\operatorname{Lie} (\tilde{H})$.
Then by $G_m' \subset \tilde{H} \subset \operatorname{SL}(V_m)$, 
we have the inclusions
$$
\frak g'_m \; \hookrightarrow \; \frak h\; \hookrightarrow \; \frak{sl}(V_m).
$$
Put $\frak k'_m:= \operatorname{Lie} (K_m')$, where
$K'_m$ is the isotropy subgroup of $K_m$ at the point $\hat{M}_m\in W_m^*$.
Then  $\frak g_m$ and $\frak g'_m$ are the complexifications of  $\frak k_m$ and
$\frak k'_m$, respectively (cf. \cite{Ca}). Put
$\frak p_m \, :=\, \sqrt{-1}\, \frak k_m$ and
$\frak p'_m \, :=\, \sqrt{-1}\, \frak k'_m$.
We further define
$$
\frak{K}_m\;  :=\,  \left \{ \bigoplus_{k=1}^{\nu_m}\,\frak{u}(V_T (\chi_k ); \rho_{h (\ell )})\right \}
\bigcap \frak{su}(V_m; \rho_{h(\ell )}),\qquad 
\frak{P}_m\, :=\, \sqrt{-1}\, \frak{K}_m.
$$
By the above inclusions of Lie algebras (see also (2.5)),
we can regard $\frak p := \sqrt{-1}\, \frak k$ as a Lie subalgebra of $\frak{P}_m$.
 Let $\omega_{\operatorname{FS}}$ be the Fubini-Study metric on $\Bbb P^* (V_m )$
defined by
$$
\omega_{\operatorname{FS}}\,:=\, (\sqrt{-1}/2\pi ) \,\partial\bar{\partial}
\log (\Sigma_{\alpha =0}^{N_m}\, |z_\alpha|^2).
$$
For each $Q \in \frak P_m$,  let $\mathcal{Q}$ be the associated
holomorphic vector field on $\Bbb P^* (V_m)$.
By the notation for $t = 0$ in Step 1 later in Proof of Assertion,
we obtain a vector field $\mathcal{Q}_{TM_m}$ on $M_m$
via the orthogonal projection  of $\mathcal{Q}$ along $M_m$ to tangential directions.
Then we have a unique real-valued function $\varphi_Q$ on 
$\Bbb P^* (V_m )$ satisfying both $\int_{\Bbb P^* (V_m )} \varphi_Q\,
\omega_{\operatorname{FS}}^{N_m} = 0$ and
$$
i^{}_{\mathcal{Q}} (\omega^{}_{\operatorname{FS}}/m)=
(\sqrt{-1}/2\pi )\,\bar{\partial}
\varphi^{}_Q.
$$ 
Let $\square^{}_{M, \operatorname{FS}} := -\,\bar{\partial}^*\bar{\partial}$
denote the Laplacian on functions on the K\"ahler manifold 
$(M, \Phi_m^*\omega_{\operatorname{FS}})$. 
Define a positive semidefinite $K_m'$-invariant inner product $(\;,\;)$ on $\frak P_m$ by
setting
\begin{align*}
&(Q_1, Q_2)\; :=\; \frac{1}{m^2}\int_{M_m} ((\mathcal{Q}^{}_1)^{}_{TM_m},
(\mathcal{Q}^{}_2)^{}_{TM_m})_{\omega^{}_{\operatorname{FS}}} \,\omega_{\operatorname{FS}}^n \\ 
&=\; \frac{\sqrt{-1}}{2\pi} \int_{M_m} \,\partial \varphi^{}_{Q_2} \wedge
\bar{\partial}\varphi^{}_{Q_1}
\wedge n\,\omega_{\operatorname{FS}}^{n-1} 
\; = \; \int_{M_m} ( \bar{\partial}\varphi^{}_{Q_1},
\bar{\partial}\varphi^{}_{Q_2})_{\omega^{}_{\operatorname{FS}}} \,
\omega_{\operatorname{FS}}^n \\
&=\; - \int_{M} \varphi^{}_{Q_1} (\square^{}_{M, \operatorname{FS}}
\varphi^{}_{Q_2})\, 
  \Phi_m^*\omega_{\operatorname{FS}}^n \, \in \Bbb R
\end{align*}
for all $Q_1$, $Q_2\in \frak P_m$. Restrict this inner product to $\frak p_m$.
Then the inner product $(\;,\; )$ on $\frak p_m$ is positive definite on $\frak p$ and hence on
$\frak p_m'$. As vector spaces, $\frak P_m$ and
$\frak p_m$ are written respectively as  orthogonal direct sums
$$
\frak P_m \, = \, \frak p \oplus \frak p^{\perp},
\qquad \; \frak p_m \, =\, \frak p_m' \oplus \frak p_m'',
$$
where $\frak p^{\perp}$ is the orthogonal complement
of $\frak p$ in $\frak P_m$, and moreover
$\frak p_m''$ is the orthogonal complements 
of
$\frak p_m'$  in $\frak p_m$ (cf. \cite{M5}).  
Hence if $Q \in \frak p^{\perp}$, then for any holomorphic vector field $\mathcal{W}$ on $M_m$,
we have
\begin{align*}
&\frac{1}{m^2}\int_{M_m} (\mathcal{Q}^{}_{TM_m}, \mathcal{W}^{}_{TM_m})_{\omega^{}_{\operatorname{FS}}}
\,\omega_{\operatorname{FS}}^n \\
&=\; \int_{M_m} (\bar{\partial}\varphi^{}_Q, \,\bar{\theta}_0 + \bar{\partial}
(\varphi_{W^1} + \sqrt{-1}\,\varphi_{W^2}))_{\omega^{}_{\operatorname{FS}}}
\,\omega_{\operatorname{FS}}^n\\
&= \; \frac{\sqrt{-1}}{2\pi} \int_{M_m}\left\{\, {\theta}_0 \wedge \bar{\partial}\varphi_Q
\,+ \, \partial (\varphi_{W^1} - \sqrt{-1}\,\varphi_{W^2}) \wedge
\bar{\partial}\varphi_Q \right\}
\wedge n\,\omega_{\operatorname{FS}}^{n-1} \;
=\; 0,
\end{align*}
where 
$i_{\mathcal{W}}  (\omega^{}_{\operatorname{FS}}/m)$ on $M_m$
 is known to be expressible as $\bar{\theta}_0 + \bar{\partial} (\varphi_{W^1}
 + \sqrt{-1}\, \varphi_{W^2})$ for some 
holomorphic $1$-form $\theta_0$ on $M_m$ and elements $W^1$, $W^2$ in $\frak p$. 
We consider the open neighbourhood (cf. \cite{M5})
$$
U_m \; :=\; \{\, X \in \frak p_m''\,;\, \zeta (\operatorname{ad} X) \frak p_m' \cap \frak p_m'' 
= \{0\}\,\}
$$
of the origin in $\frak p_m''$, where $\zeta : \Bbb R \to \Bbb R$ is a real analytic function 
defined by
$\zeta (x) := x (e^x + e^{-x})/(e^x - e^{-x})$, $x \neq 0$,  and $\zeta (0 ) =0$.
By operating $ (\sqrt{-1}/2\pi )
\,\partial\bar{\partial}
\log\,$ on both sides of (3.1), we obtain
$$
\Phi_m^*\omega_{\operatorname{FS}}\;  \equiv \; m\, \omega (\ell ),\qquad
\text{mod $q^{\ell +2}$.}
\leqno{(3.3)}
$$
For an element $X$ of $\frak P_m$ (later we further assume $X \in \frak p_m''$),
there exists a
$T$-admissible  orthonormal basis 
$\mathcal{T} := \{\tau_0, \tau_1,\dots, \tau_{N_m}\}$ for $(V_m, \rho_{h (\ell )})$
such that the infinitesimal action of $X$
on $V_m$ can
be diagonalized in the form
$$
X\cdot \tau_{\alpha} \,=\, \gamma_{\alpha}(X) \, \tau_{\alpha}
$$
for some real constants $\gamma_{\alpha} = \gamma_{\alpha}(X)$, $\alpha = 0,1,\dots, N_m$,
satisfying $\Sigma_{\alpha = 0}^{N_m} \,\gamma_{\alpha}(X) = 0$.
As in (3.2),  we consider the associated basis 
$\tilde{\mathcal{T}} = \{\tilde{\tau}_0, \tilde{\tau}_1, \dots, \tilde{\tau}_{N_m}\}$ 
for $V_m$, where $\tilde{\tau}_{k,i}  := e_{}^{-d\chi_k (F(\ell ))/2}{\tau}_{k,i}$. 
By setting
$$
\lambda_X (e^t) := \exp\, (t\, X), \qquad t \in \Bbb R,
$$
we consider the one-parameter group $\lambda_X : \Bbb R_+ \to \,
\{\, \Pi_{k=1}^{\nu_m}\operatorname{GL}(V_T (\chi_k )) \,\}\, \cap \,\operatorname{SL}(V_m)$
associated to
$X$. Then $\lambda_X (e^t) \cdot \tau_{\alpha} = 
e_{}^{t\gamma_{\alpha}}\tau^{}_{\alpha}$
 for all $\alpha$ and all $t \in \Bbb R$.
Moreover,
$$
\Phi_m^* \varphi_X \; = \; \frac{\Sigma_{\alpha =0}^{N_m} \,
\gamma_{\alpha}(X)\,  |\tilde{\tau}_{\alpha} |^2}
{m\,\Sigma_{\alpha =0}^{N_m}\,  |\tilde{\tau}_{\alpha} |^2},
\quad\qquad X \in \frak P_m.
\leqno{(3.4)}
$$
Let $\eta_m $ be the K\"ahler form on $M$ defined by
$\eta_m := (1/m)\,\Phi_m^*\omega_{\operatorname{FS}}$.
To each $X \in \frak P_m$, we can associate a real constant $c_X$ such that
$\phi_X \, :=\, c_X + \Phi_m^* \varphi_X$ 
on $M$ satisfies
$$
\int_M \phi_X \, \eta_m^{n} \; =\; 0.
$$

\medskip\noindent
{\em Proof of Assertion}:\,
 Fix an arbitrary element $0\neq X$ of $\frak p''_m$, and define
a real-valued function
$f_{X,m}(t)$ on $\Bbb R$ by
$$
f_{X,m}(t)\; := \;
\log \| \lambda_X (e^t)\cdot \hat{M}_m \|_{\operatorname{CH}(\rho_{h (\ell )})}.
$$
For this $X$, we consider the associated $\gamma_{\alpha}(X)$, $\alpha = 0,1, \dots,
N_m$,  defined in the above. 
From now on, $X$ regarded 
as a holomorphic vector field on
$\Bbb P^*(V_m)$ will be denoted by $\mathcal{X}$. 
By \cite{Z2} (see also \cite{M3}, 4.5), we have
$\ddot{f}_{X,m} (t) \geq 0$ for all $t$.
Then by \cite{M5},  Lemma 3.4, it suffices to show the existence of a real number $t^{(m)}_X$ such
that 
$$
\dot{f}_{X,m}(t^{(m)}_X ) = 0 < \ddot{f}_{X,m}(t^{(m)}_X)\quad
\text{ and }\quad t^{(m)}_{X}\cdot X \in U_m.
\leqno{(3.5)}
$$
In the below, 
real numbers $C_i$, $i = 1,2,\dots$, always mean positive real constants independent of
the choice of $m$ and $X$.
Moreover by abuse of terminology, we write $m \gg 1$,  if  $m$ satisfies $m \geq m_0$ for a sufficiently large
$m_0$ independent of the choice of $X$.
Then the proof of Assertion will be divided into 
the following eight steps:

\medskip\noindent
{\em Step $1$}.
Put $\lambda_t := \lambda_X (e^t)$ 
and $M_{m,t} := \lambda_t (M_m)$ for each $t \in \Bbb R$.
Metrically, we identify the normal bundle    of $M_{m,t}$ in
$\Bbb P^*(V_m)$  with the subbundle $TM^{\perp}_{m,t}$ of  $T\Bbb P^*(V_m)_{|M_{m,t}}$
obtained as the orthogonal complement of $TM_{m,t}$ in $T\Bbb P^*(V_m)_{|M_{m,t}}$.
Hence, $T\Bbb P^*(V_m)_{|M_{m,t}}$ is differentiably written as the direct sum
$TM_{m,t} \oplus TM^{\perp}_{m,t}$. Associated to this,
the restriction $\mathcal{X}_{|M_{m,t}}$ of $\mathcal{X}$ to $M_{m,t}$ is written as
$$
\mathcal{X}_{|M^{}_{m,t}} \; =\; \mathcal{X}_{TM^{}_{m,t}} \oplus 
\mathcal{X}_{TM_{m,t}^{\perp}}\, 
$$
for some smooth sections $\mathcal{X}_{TM^{}_{m,t}}$ and
$\mathcal{X}_{TM_{m,t}^{\perp}}$ of $TM_{m,t}$ and 
$TM_{m,t}^{\perp}$\,,
respectively. 
 Then the second derivative $\ddot{f}_{X,m}(t)$ is (see for instance \cite{M3},  \cite{PS}) 
given by
$$
\ddot{f}_{X,m}(t) \; =\; \int^{}_{M_{m,t}} 
|\mathcal{X}_{TM_{m,t}^{\perp}}|^{\,2}_{\omega^{}_{\operatorname{FS}}}
 \, \omega^n_{\operatorname{FS}}\;\, \geq \,\; 0.
\leqno{(3.6)}
$$
Since the Kodaira embedding ${\Phi}^{\tilde{\mathcal{T}}} : M \to \Bbb P^{N_m}(\Bbb C )$ 
defined by 
$$
{\Phi}^{\tilde{\mathcal{T}}} (p) := (\tilde{\tau}_0 (p) : \tilde{\tau}_1 (p): \dots :
\tilde{\tau}_{N_m}(p) )
$$
 coincides with $\Phi_m$ above up to an isometry
of $(V_m, \rho_{h (\ell )} )$, we may assume without loss of
generality that
${\Phi}^{\tilde{\mathcal{T}}}$ is chosen as $\Phi_m$.

\medskip\noindent
{\em Step $2$}.
In view of the orthogonal decomposition $\frak P_m = {\frak p}^{\perp}\oplus {\frak p}$,  
we can express $ X$ as an othogonal sum
$$
X\; =\; X' + X''
$$
for some $X' \in {\frak p}$ and $X''\in {\frak p}^{\perp}$.
Since $\omega (\ell )$ is $K$-invariant (cf.~\cite{M3}), the group $K$ acts isometrically
on 
$(V_m, \rho_{h (\ell )})$. Now, there exists a $T$-admissible orthonormal basis $\mathcal{B} :=
\{\beta_0, \beta_1,
\dots, \beta_{N_m}\}$ for $V_m$
such that the infinitesimal action of
$X''$ on $V_m$ is written as
$$
X''\cdot \beta_{\alpha} \; =\; \gamma^{}_{\alpha}(X'')\, \beta_{\alpha},
\qquad \alpha = 0,1,\dots, N_m,
$$
for some real constants $\gamma_{\alpha}^{}(X'') $, $\alpha = 0,1,\dots, N_m,$
satisfying $\Sigma_{\alpha =0}^{N_m} \,\gamma_{\alpha}(X'') = 0$.
 By the notation as in (3.2),  we consider the associated basis
$\tilde{\mathcal{B}} := \{\tilde{\beta}_0, \tilde{\beta}_1,$
\dots, $\tilde{\beta}_{N_m}\}\,$ for $V_m$. Then
$$
\phi_{X''} \; = \; \frac{\Sigma_{\alpha =0}^{N_m} \,
\hat{\gamma}_{\alpha}(X'')\, |\tilde{\beta}_{\alpha} |^2}
{m\,\Sigma_{\alpha =0}^{N_m}\,  |\tilde{\beta}_{\alpha} |^2},
\leqno{(3.7)}
$$
where $\hat{\gamma}_{\alpha}(X'') := \gamma_{\alpha}(X'')  + m \,c_{X''}$.
Now, $X'$ and $X''$ regarded 
as holomorphic vector fields on
$\Bbb P^*(V_m)$ will be denoted by $\mathcal{X}'$ and $\mathcal{X}''$,
respectively. 
Associated to the expression $T\Bbb P^*(V_m)_{|M_{m,t}}
=TM_{m,t} \oplus TM^{\perp}_{m,t}$ as differentiable vector bundles,
the restrictions $\mathcal{X}'_{|M_{m,t}}$
$\mathcal{X}''_{|M_{m,t}}$ of $\mathcal{X}'$ and $\mathcal{X}''$ 
to $M_{m,t}$ are respectively written as
$$
\mathcal{X}'_{|M^{}_{m,t}} \; =\; \mathcal{X}'_{TM^{}_{m,t}} \oplus 
\mathcal{X}'_{TM_{m,t}^{\perp}}
\quad\text{ and }\quad
\mathcal{X}''_{|M^{}_{m,t}} \; =\; \mathcal{X}''_{TM^{}_{m,t}} \oplus 
\mathcal{X}''_{TM_{m,t}^{\perp}},
$$
where $\mathcal{X}'_{TM^{}_{m,t}}$, $\mathcal{X}''_{TM^{}_{m,t}}$ are smooth sections of
$TM_{m,t}$, and $\mathcal{X}'_{TM_{m,t}^{\perp}}$
$\mathcal{X}''_{TM_{m,t}^{\perp}}$ are smooth sections of 
$TM_{m,t}^{\perp}$.
 Then by $X' \in \frak p$, we have
$$
\mathcal{X}'_{TM_{m,t}^{\perp}} \; =\; 0,
\quad\text{i.e.,}\quad  \mathcal{X}_{TM_{m,t}^{\perp}}\; =\; \mathcal{X}''_{TM_{m,t}^{\perp}}.
\leqno{(3.8)}
$$
{\em Step $3$}.
Since $T$ is irredundant, we have $\frak{g}_m'(T) + \frak t = 
\frak k^{\Bbb C}$, i.e., $\frak p_m' +  \sqrt{-1}\,\frak t_c  = {\frak p}$,
where these are equalities as Lie subalgebras of $\frak h$.
From now on until the end of this step,  as in the preceding steps, we regard both $\frak p_m'$ and
$\frak k^{\Bbb C}$ as Lie subalgebras of $\frak{sl}(V_m)$. Hence, as Lie subalgebras of
$\frak{sl}(V_m)$,  we have
$$
{\frak p}  \;  = \; \frak p_m' +  \sqrt{-1}\, \tilde{\frak t}_c .
$$
where we put $\tilde{\frak t}_c := \operatorname{Lie}(T_c)$ 
for the maximal compact subgroup $\tilde{T}_c$
of $\tilde{T} := \iota^{-1}(T)$.
Then we can write  $X' \in {\frak p}$ as a sum
$$
X' \; =\; Y + W
$$
for some $Y \in \frak p_m'$ and some $W \in  \sqrt{-1}\, \tilde{\frak t}_c$.
Note that the holomorphic vector fields $\mathcal{Y}$ and $\mathcal{W}$
on $\Bbb P^* (V_m )$ induced by $Y$ and $W$, respectively,
are tangent to $M$. 
By $[Y, W] = 0$, there exists a $T$-admissible 
orthonormal basis $\{ \sigma_0, \sigma_1, 
\dots, \sigma_{N_m}\}$ for $V_m$ such that
\begin{equation*}
\begin{cases}
&Y \cdot \sigma_{\alpha} \; =\; \gamma_{\alpha}(Y)\, \sigma_{\alpha},
\qquad \alpha = 0,1,\dots, N_m;\\
 &W \cdot \sigma_{k,i} \; =\; b_k\, \sigma_{k,i}, \qquad k = 1,2,\dots, \nu_m,
\end{cases}
\end{equation*}
for some real constants $\gamma_{\alpha} (Y)$ and $b_k$, 
where in the last equality, $\sigma_{k,i} := \sigma_{l (k,i)}$ as in (2.8).
By setting $\tilde{\sigma}_{k,i} \; =\; e_{}^{-d\chi_k (F(\ell ))/2}{\sigma}_{k,i}$, 
we later consider the basis
$\{\tilde{\sigma}_0, \tilde{\sigma}_1, \dots, \tilde{\sigma}_{N_m}\}$ for $V_m$.
Note that $\Sigma_{\alpha = 0}^{N_m} \, \gamma_{\alpha}(Y) =  \Sigma_{k = 1}^{\nu_m}\, n_k b_k = 0$.
Since both $X$ and $Y$ belong to $\frak p_m$,
it follows from $X = X' + X'' = Y + W + X''$ that
$$
\Sigma_{i=1}^{n_k}\, \gamma_{k,i}(X) \, =\, \Sigma_{i = 1}^{n_k}\, \gamma_{k,i} (Y)\, 
=\,  0\;\, \; \text{and}\;\,\;
n_k b_k \, =\, - \Sigma_{i=1}^{n_k} \,\gamma_{k, i}(X'')
\quad\text{for all $k$,}
$$
where $\gamma_{k,i} := \gamma_{l (k,i)}$ as in (2.8).
Note that  $X\in \frak p_m'' $ and $X'' \in \frak p^{\perp}$, where by $\frak p_m' 
+ \sqrt{-1}\, {\frak t} =
{\frak p}$,
the space $\frak p^{\perp}$ is perpendicular to $\frak p_m'$.
Then by $Y\in \frak p_m'$ and $X= Y + W + X''$, we have
$$
(Y, Y) + (Y, W) \, =\, (Y, X) \, =\, 0
$$
in terms of the inner product $(\;,\; )$ on $\frak P_m$. Hence
$(Y, Y)\, =\, -\, (Y, W) \, \leq \, \sqrt{(Y, Y) (W, W)}$. It now follows that
$$
\int_{M_m} | \mathcal{Y}_{|M_m} |^2_{\omega_{\operatorname{FS}}} \,\omega^n_{\operatorname{FS}}
\, =\, m^2(Y, Y) \,
\leq \,
\, m^2(W, W)  \,=\, \int_{M_m} | \mathcal{W}_{|M_m} |^2_{\omega_{\operatorname{FS}}}
\omega^n_{\operatorname{FS}}. \leqno{(3.9)}
$$ 
The integral on the right-hand side is, for $m \gg 1$,  
\begin{align*}
& \int_{M} \frac{(\Sigma_{k=1}^{\nu_m}\Sigma_{i=1}^{n_k} |\tilde{\sigma}_{k,i}|^2_{h (\ell )} )
(\Sigma_{k=1}^{\nu_m}\Sigma_{i=1}^{n_k} b_k^2 |\tilde{\sigma}_{k,i}|^2_{h (\ell )} )
- (\Sigma_{k=1}^{\nu_m}\Sigma_{i=1}^{n_k} b_k |\tilde{\sigma}_{k,i}|^2_{h (\ell )} )^2}
{(\Sigma_{k=1}^{\nu_m}\Sigma_{i=1}^{n_k} |\tilde{\sigma}_{k,i}|^2_{h (\ell )} )_{}^2 } m^n\eta_m^n \\
&\; \leq \;\; m^n \int_{M} \frac{\Sigma_{k=1}^{\nu_m}\Sigma_{i=1}^{n_k} b_k^2 |\tilde{\sigma}_{k,i}|^2_{h
(\ell )} } {\Sigma_{k=1}^{\nu_m}\Sigma_{i=1}^{n_k} |\tilde{\sigma}_{k,i}|^2_{h (\ell )}  } \,\eta_m^n\\
&\leq \;\;  \frac{n!}{2}\int_M \Sigma_{k=1}^{\nu_m}\Sigma_{i=1}^{n_k} b_k^2 |\tilde{\sigma}_{k,i}|^2_{h (\ell )} 
\,\eta_m^n 
\; \leq \;  C_1\Sigma_{k=1}^{\nu_m} n_k b_k^2
\end{align*}
for some $C_1$, where in the last two inequalities, we used Remark 2.11 in Section 2.
Hence, by setting $\hat{\gamma}_{\alpha}(Y) := \gamma_{\alpha}(Y) + m c_Y$, we see from (3.9) that
\begin{align*}
 &\int_{M} \frac{(\Sigma_{\alpha =0}^{N_m} |\tilde{\sigma}_{\alpha}|^2_{h (\ell )} )
(\Sigma_{\alpha =0}^{N_m} \hat{\gamma}_{\alpha}(Y)^2 |\tilde{\sigma}_{\alpha}|^2_{h (\ell )} )
\,-\, (\Sigma_{\alpha =0}^{N_m} \hat{\gamma}_{\alpha}(Y) |\tilde{\sigma}_{\alpha}|^2_{h (\ell )} )^2}
{(\Sigma_{\alpha =0}^{N_m} |\tilde{\sigma}_{\alpha}|^2_{h (\ell )} )_{}^2 } \, \eta_m^n 
\tag{3.10}
\\
&\;\; \; \quad \leq \; \;\;  \; q^n C_1\, \Sigma_{k=1}^{\nu_m}\, n_k b_k^2.
\end{align*}
Define real numbers $f_1$ and $f_2$ by 
\begin{equation*}
\begin{cases}
\; f_1 &:=  \int_{M} (\Sigma_{\alpha =0}^{N_m} |\tilde{\sigma}_{\alpha}|^2_{h (\ell )} )^{-1}
(\Sigma_{\alpha =0}^{N_m} \hat{\gamma}_{\alpha}(Y)^2
|\tilde{\sigma}_{\alpha}|^2_{h (\ell )})
 \,  \eta_m^n,\\
\; f_2 &:= \int_{M} \{\,(\Sigma_{\alpha =0}^{N_m} |\tilde{\sigma}_{\alpha}|^2_{h (\ell )} )^{-1}
(\Sigma_{\alpha =0}^{N_m} \hat{\gamma}_{\alpha}(Y)\,
|\tilde{\sigma}_{\alpha}|^2_{h (\ell )})\,\}^{2}
 \,  \eta_m^n.
\end{cases}
\end{equation*}
 If $f_1 \geq 2 f_2$, then by (3.10) and Remark 2.11, we have 
$$
\Sigma_{\alpha =0}^{N_m} \, \hat{\gamma}_{\alpha}(Y)^2 \; \leq\;  C_2\, \Sigma_{k=1}^{\nu_m}\,
 n_k
b_k^2, \qquad \text{if $m \gg 1$},
$$ 
for some $C_2$. Next, assume $f_1 < 2 f_2$. Then for 
$$
\phi_Y :=  (m\Sigma_{\alpha =0}^{N_m} |\tilde{\sigma}_{\alpha}|^2_{h (\ell )} )^{-1}
(\Sigma_{\alpha =0}^{N_m} \hat{\gamma}_{\alpha}(Y)
|\tilde{\sigma}_{\alpha}|^2_{h (\ell )} ), 
$$
the left-hand side of (3.10) divided by $m \gg 1$ is written as
$$
\int_{M_m} | \mathcal{Y}_{|M_m} |^2_{(\omega_{\operatorname{FS}}/m)} \, (\omega_{\operatorname{FS}}/m)^n \; =\; 
\int_M |\bar{\partial}\phi_Y |_{\eta_m}^2\, \eta_m^n\; 
\geq \; C_3  \int_M
\phi_Y^{\,2}\,\eta_m^n,
$$
for some $C_3$, because the K\"ahler manifolds $(M, \eta_m)$, $m \gg 1$, have
bounded geometry (see also Remark 2.11).
Hence, by $f_1 < 2 f_2$ and (3.10), we see that, for $m \gg 1$, 
\begin{align*}
q^n C_1\, \Sigma_{k=1}^{\nu_m}\, n_k b_k^2\; &\geq\; m \, C_3 \int_M
\phi_Y^{\,2}\,\eta_m^n \;
 =  \; C_3\, q\,  f_2 \\
& > \; C_3\, q\,  f_1/2 
\; \geq \; C_{4}\, q^{n+1} \Sigma_{\alpha =0}^{N_m} \, \hat{\gamma}_{\alpha}(Y)^2
\end{align*}
for some $C_{4}$, where in the last inequality, we used Remark 2.11.
By $\Sigma_{\alpha = 0}^{N_m}\, {\gamma}_{\alpha}(Y)$ = $0$, 
we here observe that $\Sigma_{\alpha =0}^{N_m} \, {\gamma}_{\alpha}(Y)^2\,\leq \, \Sigma_{\alpha
=0}^{N_m} \, \hat{\gamma}_{\alpha}(Y)^2$, Hence, whether $f_1 \geq 2 f_2$ or not, there 
always exists 
$C_{5}$ such that, for $m \gg 1$, 
$$
q\, \Sigma_{\alpha = 0}^{N_m} \, {\gamma}_{\alpha}(Y)^2 \; \leq \; C_{5}\, \Sigma_{k=1}^{\nu_m}\, n_k
b_k^2.
\leqno{(3.11)}
$$
{\em Step $4$}. Put $P := W + X''$. In view of \cite{Z2}, Theorem 1.6, a weighted
version of (3.4.2) in \cite{Z2} is true (cf. \cite{M3}, \cite{MW}).
Hence by $T\in \Delta_1$, we obtain
$$
\dot{f}^{}_{X,m}(0) \; = \; \dot{f}^{}_{P,m}(0) \; = \; (n+1) \, \int_M\frac{
\Sigma_{\alpha =
0}^{N_m}\, 
{\gamma}_{\alpha}(P) \,|\tilde{\beta}_{\alpha}|_{h(\ell )}^2}
{\Sigma_{\alpha =
0}^{N_m}\, |\tilde{\beta}_{\alpha}|_{h(\ell )}^2}\,
\Phi_m^* \omega_{\operatorname{FS}}^n,
\leqno{(3.12)}
$$
where ${\gamma}^{}_{k,i}(P) := b_k + {\gamma}^{}_{k,i}(X'') \, (= {\gamma}^{}_{l (k,i)}(P))$.
Let $C_{q,\ell} = 1 + \Sigma_{j =0}^{\ell} \alpha_j q^{j+1}$ 
be as in (3.1).
Then by (3.1) and (3.3), there exist a function $u_{m, \ell}$ and a $1$-form $\theta_{m, \ell }$ on $M$
such that
\begin{equation}
\begin{cases}
 &B_q (\omega (\ell ), F(\ell )) \; =\; (n!/m^n)\Sigma_{\alpha = 0}^{N_m}\,
|\tilde{\beta}_{\alpha}|_{h(\ell )}^{\,2} \; = \;
 C_{q,\ell} + u^{}_{m, \ell}\, q^{\ell  + 2};
 \\
&\eta_m =(1/m) \Phi_m^*\omega_{\operatorname{FS}} = \omega (\ell ) 
+ \theta^{}_{m,\ell}\, q^{\ell +2},
\end{cases}\tag{3.13}
\end{equation}
where we have the inequalities $\|u_{m, \ell}\|_{C^0(M)}\leq  C_{6}$ and 
$\|\theta_{m, \ell} \|_{C^0(M,
\omega_0)} \leq C_{7}$  for some $C_{6}$ and $C_{7}$ (cf. Remark 2.11; 
see also \cite{Z1}, \cite{M3}).
Hence, if $m \gg 1$,
\begin{align*}
&\;\; \;\;\;\; \frac{\,|\dot{f}_{X,m}(0)|\,}{(n+1)!}\;  \leq \; \\
& \int_M 
\frac{\left (\Sigma_{\alpha = 0}^{N_m}\, 
{\gamma}_{\alpha} (P) \, |\tilde{\beta}_{\alpha}|_{h(\ell )}^2 \right )
\left \{ 1 + \Sigma_{i = 1}^{\infty}(- u_{m, \ell}\, q^{\ell + 2}/C_{q,\ell} )^i\right \}
\left \{ \omega (\ell ) 
+ \theta_{m, \ell}\, q^{\ell +2}\right \}^n}
{C_{q,\ell}}.
\end{align*}
Here, $ \{\, 1 + \Sigma_{i = 1}^{\infty}\,(- u_{m,\ell}\, q^{\ell + 2}/C_{q,\ell} )^i \,\}
 \{\, \omega (\ell ) 
+ \theta_{m,\ell}\, q^{\ell +2}\, \}^n $ is written as $(1 + w_{m,\ell})\,\omega (\ell )^n $
for some function $w_{m, \ell}$ on $M$ such that the inequality $\| w_{m,\ell}\|_{C^0(M)} \leq C_{8}$ 
holds for some $C_{8}$.
Then by $\int_M \{\Sigma_{\alpha = 0}^{N_m}\, 
{\gamma}_{\alpha} (P)|\tilde{\beta}_{\alpha}|_{h(\ell )}^2\}\,\omega (\ell )^n =
\Sigma_{k=1}^{\nu_m}\{ e^{d\chi^{}_k(F (\ell ))}$ $\Sigma_{i=1}^{n_k}
(b_k + \gamma_{k, i }(X''))\,\}  = 0$,
we have
\begin{align*}
& |\dot{f}_{X,m}(0)| \, \leq \, (n+1)!\,q^{\ell +2}\int_M 
\left |\frac{\Sigma_{\alpha = 0}^{N_m}\, 
{\gamma}_{\alpha}(P)\, |\tilde{\beta}_{\alpha}|_{h(\ell )}^2}{C_{q,\ell}}\right |\cdot
 |w_{m, \ell} | \;\omega (\ell )^n\\
& =  (n+1)!q^{\ell +2-n}\int_M \{1 + (u_{m, \ell} q^{\ell +2}/C_{q,\ell} ) \}
\left |\frac{ \Sigma_{\alpha = 0}^{N_m} 
{\gamma}_{\alpha}(P) |\tilde{\beta}_{\alpha}|_{h(\ell )}^2}{\Sigma_{\alpha = 0}^{N_m}
|\tilde{\beta}_{\alpha}|_{h(\ell )}^{\,2} }\right | \cdot
 \left |\frac{w_{m, \ell}}{n!}  \right | \omega (\ell )^n.
\end{align*}
In view of (3.4), by setting $\hat{\phi} := (m\,\Sigma_{\alpha = 0}^{N_m}\,
|\tilde{\beta}_{\alpha}|^{\,2} )^{-1}(\Sigma_{\alpha = 0}^{N_m}\, 
{\gamma}_{\alpha}(P)\, |\tilde{\beta}_{\alpha}|^2)$,
there exist  $C_{9}$ and
$C_{10}$ 
such that,  for $m \gg 1$, 
$$
|\dot{f}_{X,m}(0)|\;\leq \; q^{\ell +1-n} C_{9}\| \hat{\phi} \|_{L^1(M,\omega
(\ell ))}
\;\leq\; q^{\ell +1-n} C_{10}\| \hat{\phi} \|_{L^2(M,\omega  (\ell
))}. \leqno{(3.14)}
$$

\medskip\noindent
{\em Step $5$}.
 Note that
$0 \leq  n_k^2 b_k^2 = \{\,\Sigma^{n_k}_{i=1} \gamma_{k,i}(X'')\,\}^2 
\leq  n_k \Sigma_{i=1}^{n_k} \gamma_{k,i}(X'')^2$ holds for all $k$ by the 
Cauchy-Schwarz inequality. 
Hence
$$
\Sigma_{k=1}^{\nu_m}\,n_k b_k^2 
\; \leq \; \Sigma_{\alpha = 0}^{N_m}\, \gamma_{\alpha}(X'')^2.
\leqno{(3.15)}
$$
From $\Sigma_{\alpha = 0}^{N_m}\,
\gamma_{\alpha}(X'') = 0$ and $\hat{\gamma}_{\alpha}(X'')  = \gamma_{\alpha}(X'') +  m\, c_{X''}$, 
it follows that $\Sigma_{\alpha = 0}^{N_m} \,\hat{\gamma}_{\alpha}(X'')^2 \, =\, 
(N_m + 1) (m \,   c_{X''})^2 \, +\, \Sigma_{\alpha = 0}^{N_m} \, {\gamma}_{\alpha}(X'')^2$. In particular,
$$
\Sigma_{\alpha = 0}^{N_m} \, {\gamma}_{\alpha}(X'')^2 \; \leq \; 
\Sigma_{\alpha = 0}^{N_m} \,\hat{\gamma}_{\alpha}(X'')^2.
\leqno{(3.16)}
$$
Since ${\gamma}_{k,i}(P) = b_k + \gamma_{k,i}(X'')$, 
(3.15) and (3.16) above imply that
$$
\Sigma_{\alpha =0}^{N_m}\, {\gamma}_{\alpha}(P)^2 \; \leq \; 
2 \{\,\Sigma_{k=1}^{\nu_m}\,n_k b_k^2 \, +\, \Sigma_{\alpha =
0}^{N_m}\, \gamma_{\alpha}(X'')^2\,\}\; \leq \;   4\Sigma_{\alpha = 0}^{N_m} \,
\hat{\gamma}_{\alpha}(X'')^2.
\leqno{(3.17)}
$$
By $X = Y + P$, we have $\Sigma_{\alpha =0}^{N_m} \, \gamma_{\alpha}(X)^2 \,
\leq  \, 2 \{\,\Sigma_{\alpha =0}^{N_m} \, \gamma_{\alpha}(Y)^2\, +\, \Sigma_{\alpha =0}^{N_m} \, 
{\gamma}_{\alpha}(P)^2\,\}$, because $\gamma_{\alpha}(X) = \gamma_{\alpha}(Y) +
\gamma_{\alpha}(P)$. Hence, by (3.11), (3.15), (3.16) and (3.17), we obtain 
\begin{align*}
&q\Sigma_{\alpha =0}^{N_m} \, \gamma_{\alpha}(X)^2\;  \leq \; 
2q\Sigma_{\alpha =0}^{N_m} \, \gamma_{\alpha}(Y)^2\, +\, 2q\Sigma_{\alpha =0}^{N_m}
\, {\gamma}_{\alpha}(P)^2          \tag{3.18}       \\
&\leq \, 2 C_5 \Sigma_{k =1}^{\nu_m} \, n_k\, b_k^2 \, +\, 
\Sigma_{\alpha =0}^{N_m}
\, {\gamma}_{\alpha}(P)^2        
\; \leq \; C_{11} \,\Sigma_{\alpha =0}^{N_m}
\, \hat{\gamma}_{\alpha}(X'')^2
\end{align*}
for $m \gg 1$, where we put $C_{11} := 4 + 2 C_5$.
Fix a positive real number $\ell_0 $  independent of the choice of $m$ and $X$. Put
$\delta_0:= q^{1/2 + \ell_0} /\sqrt{\Sigma_{\alpha =0}^{N_m}\, \hat{\gamma}_{\alpha}(X'')^2}$.
Then by (3.18), we have $0 < \delta_0 < \sqrt{C_{11}}\,/\bar{\gamma}$,
where $\bar{\gamma} := \max \{ \, |\gamma_{\alpha}(X)| \,;\, \alpha  = 0,1,\dots, N_m\,\}$.
In view of  Step 1 of \cite{M5}, Section 4, 
 by assuming $|t| \leq \delta_0$, we see that the family of K\"ahler manifolds 
$(M, q\,\Phi_m^* \lambda^*_t \omega_{\operatorname{FS}} )$ have bounded geometry.

\medskip\noindent
{\em Step $6$}. At the beginning of this step, 
we shall show the inequality (3.19) below as an analogue of \cite{PS}, (5.9), 
by proving that an argument of Phong and
Sturm \cite{PS} for
$\dim H = 0$ is valid also for $\dim H > 0$. 
To see this, we consider the following exact sequence of holomorphic vector bundles
$$
0 \to TM_{m,t} \to T\Bbb P^*(V_m)_{|M_{m,t}} \to TM_{m,t}^{\perp} \to 0,
$$
where $TM_{m,t}^{\perp}$ is regarded as the normal bundle of $M_{m,t}$ 
in $\Bbb P^*(V)$. 
The pointwise estimate (cf. \cite{PS}, (5.16)) of the second fundamental form for this exact sequence
has nothing to do with $\dim H$, and as in \cite{PS}, (5.15), it gives the inequality
$$
\int_{M_{m,t}} |\mathcal{X}''_{TM_{m,t}^{\perp}}
 |^2_{\omega_{\operatorname{FS}}} \,\omega^n_{\operatorname{FS}}
\; \geq \; C_{12} 
\int_{M_{m,t}} |\,\bar{\partial}
\mathcal{X}''_{TM_{m,t}^{\perp}} |^2_{\omega_{\operatorname{FS}}}
\,\omega^n_{\operatorname{FS}}.
$$
for some $C_{12}$.
Let $\mathcal{A}_{}^{0,p} (T_M)$, $p$ = 0,1, denote the sheaf of germs of smooth
$(0,p)$-forms on $M$ with values in the holomorphic tangent bundle $TM$ of $M$,
and endow $M$ with the K\"ahler metric $(1/m) \Phi_m^* \lambda_t^*
\omega_{\operatorname{FS}}$. 
We then consider the operator 
$\square^{}_{TM} := -\bar{\partial}_{}^{\#}\,\bar{\partial}$ on
$\mathcal{A}^{0,0} (T_M) $, where
$\bar{\partial}_{}^{\#}:\mathcal{A}^{0,1} (T_M) 
\to \mathcal{A}^{0,0} (T_M) $ is the formal adjoint of $\bar{\partial} : \mathcal{A}^{0,0} (T_M) 
\to \mathcal{A}^{0,1} (T_M) $.
Since by Step 1, the K\"ahler metrics $q \Phi_m^* \lambda_t^*
\omega_{\operatorname{FS}}$   has bounded geometry,
the first positive eigenvalue $\lambda_1$ of
the operator 
$-\square^{}_{TM}$ on
$\mathcal{A}^{0,0} (T_M) $ is bounded from below by $C_{13}$. 
Hence, by $X'' \in \frak p^{\perp}$,
$$
\int_{M_{m,t}}  |\bar{\partial}
\mathcal{X}''_{TM_{m,t}} |^2_{(\omega_{\operatorname{FS}}/m)}
\, (\omega_{\operatorname{FS}}/m )_{}^n \; \geq \; 
\lambda_1 \int_{M_{m,t}}|\mathcal{X}''_{TM^{}_{m,t}}|^{\,2}_{(\omega^{}_{\operatorname{FS}}/m)}
\, (\omega_{\operatorname{FS}}/m)_{}^n.
$$
Since $\bar{\partial}\mathcal{X}''_{TM_{m,t}^{\perp}} = - 
\bar{\partial}\mathcal{X}''_{TM_{m,t}}$,
by $\lambda_1 \geq C_{13}$, it now follows
that
\begin{align*}
\ddot{f}_{X,m}(t) \;\; &=\;\;\int_{M_{m,t}}
|\mathcal{X}''_{TM_{m,t}^{\perp}}|^{\,2}_{\omega^{}_{\operatorname{FS}}}
\omega^{n}_{\operatorname{FS}} \tag{3.19} \\
&\geq\;\;
C^{}_{12}C^{}_{13}\,q\int_{M_{m,t}}|\mathcal{X}''_{TM^{}_{m,t}}|^{\,2}_{\omega^{}_{\operatorname{FS}}}
\omega^{n}_{\operatorname{FS}}.
\end{align*}
In view of the equality $|\mathcal{X}''_{TM^{}_{m,t}}|^{\,2}_{\omega^{}_{\operatorname{FS}}}
+ |\mathcal{X}''_{TM_{m,t}^{\perp}}|^{\,2}_{\omega^{}_{\operatorname{FS}}}
= |\mathcal{X}''_{|M^{}_{m,t}}|^{\,2}_{\omega^{}_{\operatorname{FS}}}$,
by adding the integral
$C_{12}C_{13}\, q \int_{M_{m,t}} |\mathcal{X}''_{TM^{\perp}_{m,t}} |^2_{\omega_{\operatorname{FS}}}
\omega_{\operatorname{FS}}^n$ to both sides of (3.19) and by dividing the resulting inequality
 by
$(1+ C_{12}C_{13} q)$, we see that, for some $C_{14}$ and $C_{15}$, 
\begin{align*}
\ddot{f}_{X,m} (t) \; &= \; \int_M \Phi_m^*\lambda_t^*\left (\,
|\mathcal{X}''_{TM_{m,t}^{\perp}}|^{\,2}_{\omega^{}_{\operatorname{FS}}}
\,\omega_{\operatorname{FS}}^n\right )   \tag{3.20}  \\
&\geq \; C_{14} \, q\, \int_M \Phi_m^*\lambda_t^*\left (\,
|\mathcal{X}''_{|M_{m,t}}|^{\,2}_{\omega^{}_{\operatorname{FS}}}
\,\omega_{\operatorname{FS}}^n\right )    \\
& \geq  \; C_{15}\, q\,
\int_{M_m}
|\mathcal{X}''_{|M_{m}} |_{\omega^{}_{\operatorname{FS}}}^{\,2}\,\omega_{\operatorname{FS}}^n
\; \geq\;  C_{15} \, q\, 
\int_M \, \Theta\, \Phi_m^* \omega_{\operatorname{FS}}^n,
\end{align*}
where $\Theta := (\Sigma_{\alpha =0}^{N_m} |\tilde{\beta}_{\alpha}|{}^2 )_{}^{-2}
\{\, (\Sigma_{\alpha =0}^{N_m}\, |\tilde{\beta}_{\alpha}|{}^2)
(\Sigma_{\alpha =0}^{N_m}\,\hat{\gamma}_{\alpha}(X'')^2|\tilde{\beta}{}_{\alpha}|^2)\, -\, 
(\Sigma_{\alpha =0}^{N_m}\,\hat{\gamma}_{\alpha}(X'')$ $|\tilde{\beta}_{\alpha}|{}^2)^2\,\}\,$
is nonnegative everywhere on $M$.
Then by (3.14) and (3.20), 
\begin{equation}
\begin{cases}
\; \dot{f}_{X,m}(\delta_0 )  &\geq\;  \dot{f}_{X,m}(0) + C_{15}\,\delta_0\, q \,\int_M
\Theta 
\Phi_m^* \omega_{\operatorname{FS}}^n  \\ 
&\geq\;
-q^{\ell +1 -n}C_{10}\| \hat{\phi} \|_{L^2(M, \omega (\ell ) )}\,+\, 
C_{15}\, \delta_0 \, q \,\int_M  \Theta 
\Phi_m^* \omega_{\operatorname{FS}}^n,\\
 \dot{f}_{X,m}(-\delta_0) &\leq\;  \dot{f}_{X,m}(0) -  C_{15}\, \delta_0\, q \,\int_M \Theta 
\Phi_m^* \omega_{\operatorname{FS}}^n  \\ 
&\leq\;  q^{\ell +1 -n}C_{10}\|\hat{\phi} \|_{L^2(M, \omega (\ell ) )}\,- \,
C_{15}\, \delta_0\, q \,\int_M  \Theta 
\Phi_m^* \omega_{\operatorname{FS}}^n.
\end{cases}\tag{3.21}
\end{equation}
By (3.20) and \cite{M5}, Lemma 3.4,
the proof of Main Theorem
is reduced to showing the following three conditions for all $m \gg 1$:
$$
\text{i) $\dot{f}_{X,m}(\delta_0 ) >0 > \dot{f}_{X,m}(-\delta_0)$, \quad  ii) $\int_M
\Theta\Phi_m^*\omega_{\operatorname{FS}}^n >0$, \quad  iii)  $t^{(m)}_X \cdot X \in U_m$.}
$$
Since  iii) follows from Remark 3.31 below, we have only to prove i) and ii).
Then by (3.21), it suffices to show the following for all $m \gg 1$:
$$
C_{15}\,\delta_0\,q\,\int_M  \Theta
\Phi_m^* \omega_{\operatorname{FS}}^n \,  - \, C_{10}\,
q^{\ell +1 -n}\|\hat{\phi} \|_{L^2(M,
\omega (\ell ) )}\; >\; 0.
\leqno{(3.22)}
$$
Let us define real numbers $\hat{e}_1$, $\hat{e}_2$, ${e}_1$, $e_2$ by setting
\begin{align*}
\hat{e}_1 &:= \int_M
\frac{\Sigma_{\alpha
=0}^{N_m} \hat{\gamma}_{\alpha}(X'')^{2}|\tilde{\beta}_{\alpha}|^2}
{\Sigma_{\alpha =0}^{N_m} |\tilde{\beta}_{\alpha}|^2}\,\omega (\ell )^n,\; \;
\hat{e}_2 := \int_M  \left (
\frac{\Sigma_{\alpha =0}^{N_m} \hat{\gamma}_{\alpha}(X'') |\tilde{\beta}_{\alpha}|^2}
{\Sigma_{\alpha =0}^{N_m} |\tilde{\beta}_{\alpha}|^2}
\right )^2\omega (\ell )^n, \\
{e}_1 &:= \int_M
\frac{\Sigma_{\alpha
=0}^{N_m}{\gamma}_{\alpha}(P)^{2}|\tilde{\beta}_{\alpha}|^2}
{\Sigma_{\alpha =0}^{N_m} |\tilde{\beta}_{\alpha}|^2}\omega (\ell )^n,\; \quad
e_2 := \int_M  \left (
\frac{\Sigma_{\alpha =0}^{N_m} {\gamma}_{\alpha}(P) |\tilde{\beta}_{\alpha}|^2}
{\Sigma_{\alpha =0}^{N_m} |\tilde{\beta}_{\alpha}|^2}
\right )^2\omega (\ell )^n.
\end{align*}
By the Cauchy-Schwarz inequality, we always have $\hat{e}_1 \geq \hat{e}_2$ and ${e}_1 \geq e_2$.
Now, the following cases are possible:

\smallskip\noindent
Case 1: \;\;  $\hat{e}_1 > 2\, \hat{e}_2$, \qquad\qquad Case 2:\;\;  $ \hat{e}_1 \leq 2\, \hat{e}_2$.

\smallskip\noindent
In view of the identities in (3.13), we can write
\begin{align*}
\hat{e}_1 \quad
&= \; q^n n!\int_M\,\frac{ \Sigma_{\alpha =0}^{N_m}\,\hat{\gamma}_{\alpha}(X'')^{2}
|\tilde{\beta}_{\alpha}|_{h(\ell
)}^2 }
{ 1 + \Sigma_{\alpha =0}^{\ell}\, \alpha_k q^{k+1}
+ u_{m, \ell}\, q^{\ell +2} }\;
\omega (\ell )^n
,\\
{e}_1 \quad
&= \; \,q^n n!\int_M\,\frac{\Sigma_{\alpha =0}^{N_m}\,{\gamma}_{\alpha}(P)^{2}
|\tilde{\beta}_{\alpha}|_{h(\ell
)}^2}
{ 1 + \Sigma_{\alpha =0}^{\ell}\, \alpha_k q^{k+1} + u_{m, \ell}\, q^{\ell +2}}
\;\omega (\ell )^n,\\ 
\int_M \Theta \Phi_m^* \omega_{\operatorname{FS}}^n &= \;
m^n \int_M \, \Theta \{\omega (\ell ) + \theta_{m, \ell}\, q^{\ell + 2}\}^n,
\end{align*}
and hence, given a positive real number $0< \varepsilon \ll 1$,  both  $\hat{e}_1$ 
and $\int_M \Theta\Phi_m^*\omega_{\operatorname{FS}}^n$ above 
are estimated,
for all $m \gg 1$,  by
\begin{align*}
&(1- \varepsilon )\, q^n_{}\,  \{\Sigma_{\alpha
=0}^{N_m}\, \hat{\gamma}_{\alpha}(X'')^2\} \; \leq \; \hat{e}_1/n! \; \leq 
(1+ \varepsilon ) \, q^n_{}\,  \{\Sigma_{\alpha
=0}^{N_m}\, \hat{\gamma}_{\alpha}(X'')^2\},
 \tag{3.23}  \\
&(1- \varepsilon )\, q^n_{}\,  \{\Sigma_{\alpha
=0}^{N_m}\, {\gamma}_{\alpha}(P)^2\} \; \leq \; {e}_1/n! \; \leq 
(1+ \varepsilon ) \, q^n_{}\, \{\Sigma_{\alpha
=0}^{N_m}\, {\gamma}_{\alpha}(P)^2\},
 \tag{3.24}  \\
&(1- \varepsilon )\, q^{-n}_{} \int_M \Theta \omega (\ell )^n\; \leq \;\int_M \Theta \Phi_m^*
\omega_{\operatorname{FS}}^n\; \leq\; (1+ \varepsilon ) \,q^{-n}_{} \int_M \Theta \omega (\ell )^n,
\tag{3.25}
\end{align*}
where  we used Remark 2.11.
Moreover,  we can write $e_2$ in the form
$$
q^{-2} \|\hat{\phi}\|^{\,2}_{L^2(M, \omega (\ell )) }
\; =\;  e_2 \; \leq \; {e}_1.
\leqno{(3.26)}
$$
{\em Step $7$}. We first consider Case 1. 
Then from (3.17), (3.23), (3.24), (3.25), (3.26), $\hat{e}_2 <
\hat{e}_1/2$ and the definition of $\delta_0$, it follows that
\begin{align*}
&\text{L.H.S. of (3.22)}\; \geq \; (1-\varepsilon)\, C_{15}\,\delta_0\,q^{1-n} \,\int_M  \Theta 
\omega (\ell )^n \,  - \, q^{\ell +1  -n}C_{10}\|\hat{\phi} \|_{L^2(M,
\omega (\ell ) )}\; \\
&\geq \; (1-\varepsilon )\, C_{15}\,q^{1-n} \,\delta_0\,(\hat{e}_1 - \hat{e}_2 )\,  - \, q^{\ell +2
-n} C_{10}\sqrt{e_1} \\
&\geq \; (1-\varepsilon )\, C_{15}\,\delta_0\,q^{1-n} \, \hat{e}_1/2 
- \, (1 + \varepsilon )^{1/2} C_{10} \,q^{\ell +2
-\frac{n}{2}} \{\, n! \Sigma_{\alpha = 0}^{N_m} \gamma_{\alpha} (P)^2 \,\}^{1/2}\\
&\geq\; (1-\varepsilon )^2C_{15}\,\delta_0\, q \, \{\,\Sigma_{\alpha
=0}^{N_m}\, \hat{\gamma}_{\alpha}(X'')^2\,\} \, n!/2
\\ 
&\qquad \qquad -\,  2(1+ \varepsilon )^{1/2}C_{10}\,q^{\ell +2
-\frac{n}{2}} \{\,  n! \,\Sigma_{\alpha
=0}^{N_m}\, \hat{\gamma}_{\alpha}(X'')^2\,\}^{1/2} \\
&\geq \,\{n!\,\Sigma_{\alpha
=0}^{N_m} \hat{\gamma}_{\alpha}(X'')^2\,\}^{1/2}
\left \{ \sqrt{n!}\, (1-\varepsilon )^2C_{15} q^{\ell_0 +\frac{3}{2}}/2
 - 2(1+ \varepsilon )^{1/2}C_{10}q^{\ell +2
-\frac{n}{2}} \right \},
\end{align*}
for $m \gg 1$.
Now we see that, if $\ell >
(n-1)/2 + \ell_0 $, then $q^{\ell + 2 - \frac{n}{2}}/q^{\ell_0+\frac{3}{2}}$
converges to $0$ as $m \to \infty$.
Thus, 
if  $m \gg 1$,  then by choosing $\ell$ such that
$\ell >
(n-1)/2 + \ell_0 $, we now see from the computation above that L.H.S. of (3.22) is positive, 
as required.

\medskip\noindent
{\em Step $8$}. Let us finally consider Case 2. 
For each fixed $\ell$, the K\"ahler form $\eta_m$ converges to $\omega_0$, as
$m\to\infty$, in
$C^j(M)$-norm for all positive integers 
$j$
 (cf. \cite{Z1}, \cite{M3}; see also Remark 2.11).
Note that $ X''\in {\frak{p}}^{\perp}$. In view of $\int_M \phi_{X''} \eta_m^{\,n} = 0$, we see
that
$$
\|\phi_{X''}\|^2_{L^2(M, \eta_m )}\,
\leq \,  C_{16} \, \|\bar{\partial} \phi_{X''} \|^2_{L^2(M,\eta_m )}\, 
=\,C_{16}\,\| \Phi_m^*\mathcal{X}''_{TM_m} \|^2_{L^2(M,\eta_m )},
\leqno{(3.27)}
$$
for some $C_{16}$, where by abuse of terminology, the differential
$(\Phi_m^{-1})_* : TM_m \to TM$ is denoted by $\Phi_m^*$. Moreover, by (3.19) applied to $t = 0$, 
we obtain 
$$
 \|\,\Phi_m^*\mathcal{X}''_{TM_{m}}
\|^{\,2}_{L^2(M, \eta_m )}
\,\leq\, (C_{12}C_{13})^{-1}q^{-1}\|\Phi_m^*\mathcal{X}''_{TM_{m}^{\perp}}\|^{\,2}_{L^2(M,
\eta_m )}.\leqno{(3.28)}
$$
From now on until the end of this proof, we assume that $m \gg 1$.
By (3.27) together with (3.28) and (3.7), there exist $C_{17}$ and $C_{18}$ such that
\begin{equation}
\begin{cases}
&\| \Phi_m^*\mathcal{X}''_{TM_m^{\perp}}\|_{L^2(M, \eta_m )}\;\geq \;  C_{17}\,q^{1/2}
\,\| \phi_{X''} \|_{L^2(M, \eta_m )}\\
& \geq \;  C_{18}\,q^{1/2}
\,\| \phi_{X''} \|_{L^2(M, \omega (\ell ) )}\;
= \; C_{18}\,q^{3/2}\, \sqrt{\hat{e}_2}.
\end{cases}\tag{3.29}
\end{equation}
We now observe the pointwise estimate
$q^{1/2}|\mathcal{X}''_{|M_m}|_{\omega_{\operatorname{FS}}}$ 
= $|\mathcal{X}''_{|M_m}|_{\eta_m}$
$\geq |\mathcal{X}''_{TM_m^{\perp}}|_{\eta_m}$.
Hence by (3.20) and (3.29), we obtain
\begin{equation}
\begin{cases}
&\ddot{f}_{X,m}(t) \; \geq \; C_{15} \int_{M_m}\left ( q^{1/2}
| \mathcal{X}''_{|M_m} |_{\omega_{\operatorname{FS}}}\right )^2 
\omega^n_{\operatorname{FS}}\;\\
&\geq\;
C_{15}\, q^{-n}
\| \Phi_m^*\mathcal{X}''_{TM_m^{\perp}}\|^2_{L^2(M, \eta_m )} 
\;\geq\; C_{19}\, q^{3-n} \hat{e}_2,
\end{cases}\tag{3.30}
\end{equation}
for some $C_{19}$.
As in deducing (3.21) from (3.14) and (3.20), we obtain by
 (3.14) and (3.30) the inequalities
$$
\dot{f}_{X,m}(\delta_0 ) \;\geq\; R
\quad \text{ and }\quad 
\dot{f}_{X,m}(-\delta_0 ) \;\leq\; -R,
$$
where $R:= - q^{\ell + 1 - n}C_{10} \| \hat{\phi} \|_{L^2 (M, \omega (\ell ) )}
+  C_{19}\,\delta_0\,q^{3-n}\hat{e}_2$.
Hence, it suffices to show that $R > 0$. In view of the definition of $\hat{\phi}$ and $e_2$,
we see from $\hat{e}_1 \leq 2 \hat{e}_2$ and (3.26) that
$$
R\; =\; C_{19}\,\delta_0\,q^{3-n} \hat{e}_2\,-\, C_{10}\,q^{\ell + 2-n}\sqrt{e_2} \; \geq \; 
C_{19} \,\delta_0\,q^{3-n}\hat{e}_1/2\, -\, C_{10}\, q^{\ell + 2-n}\sqrt{{e}_1}. 
$$
Here by (3.23) and (3.24), we obtain 
\begin{align*}
&\frac{\delta_0 q^{3-n}\hat{e}_1}{q^{\ell + 2-n}\sqrt{{e}_1}} \;
=\; q^{3/2 + \ell_0 - \ell}\, \sqrt{\frac{\hat{e}_1}{{e}_1}}\sqrt{ 
\frac{\hat{e}_1}{\Sigma_{\alpha =0}^{N_m}\hat{\gamma}_{\alpha}(X'')^2} } \\
&\geq \; C_{20}\, q^{(3+n)/2 + \ell_0 - \ell}
\sqrt{\frac{\Sigma_{\alpha =0}^{N_m}\,\hat{\gamma}_{\alpha}(X'')^2}{\Sigma_{\alpha
=0}^{N_m}\, {\gamma}_{\alpha}(P)^2}}
\; \geq \; \frac{C_{20}}{2} \, q^{(3+n)/2 + \ell_0 - \ell}
\end{align*}
for some $C_{20}$, where the last inequality follows from (3.17).
Therefore, by choosing $\ell$ such that $\ell > (3+ n)/2 + \ell_0$, we now conclude that
$R >0$ for $m \gg 1$, as required. 
\qed

\medskip\noindent
{\em Remark $3.31$}. In the above proof, it is easy to check the condition  iii) in  Step 6 as follows:
In view of $|t_X^{(m)}| <  \delta_0 $,  it suffices to show that, if $m \gg 1$, then
$$
t\cdot X \in  U_m, 
\qquad \text{for all $(t, X) \in \Bbb R \times \frak p_m''$ with $|t| < \delta_0$.}
\leqno{(3.32)}
$$
For each $Q \in \frak p$, let $u_Q\in C^{\infty}(M)_{\Bbb R}$ denote the Hamiltonian function for the holomorphic vector field 
$Q$  on
the K\"ahler manifold
$(M, \omega_0 )$ characterized by the equalities
$$
i^{}_Q \omega^{}_0 \; =\; (\sqrt{-1}/2\pi )\,\bar{\partial}u^{}_Q
\quad \text{ and } \quad \int_M u^{}_Q \omega_0^{\,n} \; =\; 0.
$$
Define compact subsets $\Sigma (\frak p_m' )$,  $\Sigma (\frak p )$ 
of $\frak p$ by setting  
\begin{equation*}
\begin{cases}
\quad \Sigma (\frak p_m' ) &:= \{\, Q \in \frak p'_m\,;\,
\|\bar{\partial} u_Q\|^{}_{L^2(M, \omega_0 )}
 = 1\,\}, \\
\quad \; \Sigma (\frak p ) &:= \{\, Q \in \frak p\,;\,
\|\bar{\partial} u_Q\|^{}_{L^2(M, \omega_0 )}
 = 1\,\}.
\end{cases} 
\end{equation*}
Choose an orthonormal basis 
$\mathcal{S}:= \{s_0, s_1, \dots, s_{N_m}\}$
for the Hermitian
vector space $(V_m, \rho_{h (\ell )})$.
For the space $\mathcal{H}_m$ of all Hermitian matrices of order $N_m + 1$, define a norm
$$
\mathcal{H}_m \to \Bbb R_{\geq 0}, \qquad A = (a_{\alpha\beta} ) \,\mapsto\, 
\| A\|^{}_m := \sqrt{\operatorname{tr} A^* A} = \sqrt{\Sigma_{\alpha, \beta}\, |a_{\alpha\beta}|^2}
$$
on $\mathcal{H}_m$. Let $m \gg 1$.
The infinitesimal action of $\frak p_m$ on $V_m$ is given by
$$
Q \cdot s_{\beta} \,=\, \Sigma_{\alpha =0}^{N_m}\;s_{\alpha} \, \gamma_{\alpha\beta}(Q),
\qquad\quad  Q \in \frak p_m, 
$$
where $\gamma_Q =(\gamma_{\alpha\beta}(Q)) \in \mathcal{H}_m$ 
denotes the representation  matrix of $Q$ on $V_m$ with respect to $\mathcal{S}$. 
Let $X \in \frak p_m''$, and let $\delta_0$ be as in Step 5 above.
For $t\in \Bbb R$ with $|t| < \delta_0$, we put $\tilde{X}:=\, t\,X$.
In order to prove (3.32) above, it suffices to show
$$
\zeta (\operatorname{ad} \tilde{X}) Q \notin \frak p_m''
\qquad \text{ for all } Q \in \Sigma (\frak p_m' ).  
\leqno{(3.33)}
$$
Let $Q\in\Sigma (\frak p_m' )$. For a suitable choice of a
basis  $\mathcal{S}$ as above, we may assume
that the representation matrix $\gamma_Q$ of $Q$  is a real diagonal matrix.  Note also that $\operatorname{tr} \gamma_Q = 0$.
Let $\Phi_m : M \to \Bbb P^{N_m}(\Bbb C )$ be the Kodaira embedding of $M$ defined by (cf. (3.2)\,)
$$
\Phi_m (p) := (\tilde{s}_0(p): \tilde{s}_1(p): \dots : \tilde{s}_{N_m}(p)).
$$
In view of the definition $\eta_m :=\Phi_m^*\omega_{\operatorname{FS}}/m$ of $\eta_m$,
the Hamiltonian function  $\phi_Q$  
on $(M,\eta_m )$ associated to the holomorphic vector field $Q$ is expressed in the form
$$
\phi_Q\,  =\,
 (\Sigma_{\alpha = 0}^{N_m} \, \hat{\gamma}_{\alpha\alpha} (Q) \,|\tilde{s}_{\alpha}|^2)
/(m\Sigma_{\alpha = 0}^{N_m}\, |\tilde{s}_{\alpha}|^2).
$$
We define $\hat{\gamma}_Q := (\hat{\gamma}_{\alpha\beta}(Q))\in \mathcal{H}_m$ by
setting
 $\hat{\gamma}_{\alpha\beta}(Q) := \{{\gamma}_{\alpha\alpha}(Q)\, + \,m \,c^{}_Q\}
\delta^{}_{\alpha
\beta}$  for Kronecker's delta $\delta_{\alpha\beta}$. As in deducing (3.16) from 
$\hat{\gamma}_{\alpha}(X'') = \gamma_{\alpha}(X'') \,+\,
m \,c_{X''}$, we easily see that
$$
\| \gamma_Q \|^2_m \; \leq \; \| \hat{\gamma}_Q \|^2_m.
$$
Recall that $\eta_m$ is expressible as $\omega_0 + (\sqrt{-1}/2\pi )\, q \,\partial\bar{\partial}\xi_m$ 
for some real-valued
smooth function 
$\xi_m$ on $M$ such that 
$$
\| \xi_m \|_{C^3 (M)} \leq C_{21},
\leqno{(3.34)}
$$
where all $C_j$'s in this remark are positive constants independent of the choice of $m$, $X$ and $Q$. 
We now observe that 
$$
\phi_Q \; =\; u_Q\,  +\,  q\,(Q \xi_m ). \leqno{(3.35)}
$$
Note that $Q \in \Sigma (\frak p_m') \subset \Sigma (\frak p )$.
Since $Q$ sits in the compact set $\Sigma (\frak p )$, 
and since $\Sigma (\frak p )$ is independent of the choice of $m$,  there exist $C_{22}$ and $C_{23}$
such that 
$$
0\; <\; C_{22} \; \leq  \; 
\int_M u_Q^{\,2}\omega_0^n \; \left (= \int_M \phi^{\,2}_Q \eta_m^{\,n} \right ) \leq \; C_{23}.
$$
Note that
both $\eta_m$ and $\omega (\ell )$ converge to $\omega_0$ as
$m \to \infty$ (see (3.3) and the statement at the beginning of Step 8).
Note also that, by Remark 2.11, the function $(n!/m^n)\Sigma_{\alpha = 0}^{N_m}\, |\tilde{s}_{\alpha}|_{h (\ell )}^2$ 
on $M$ converges uniformly to $1$, as $m \to \infty$. Again by Remark 2.11,
it now follows from the Cauchy-Schwarz inequality  
that, for $m \gg   1$,
\begin{align*}
&\| \gamma_Q \|_m^{\,2} \; =\; \sum_{\alpha = 0}^{N_m}\, \gamma_{\alpha\alpha}(Q)^2\; 
\geq \; C_{24}\,
m^n\int_M \frac{\Sigma_{\alpha = 0}^{N_m} \,\gamma_{\alpha\alpha}(Q)^2\,
|\tilde{s}_{\alpha}|^2}{\Sigma_{\alpha =
0}^{N_m}\, |\tilde{s}_{\alpha}|^2}
\, \omega (\ell )^n \;  \tag{3.36} \\
&\geq  \; C_{24}\, 
m^{n+2}\int_M\, (\Phi_m^*\varphi_Q)^2 \,\omega (\ell )^n \;  \geq \; 
C_{25}\, m^{n+2}\int_M\, (\Phi_m^*\varphi_Q)^2\, \eta_m^n
\end{align*}
for some $C_{24}$ and $C_{25}$, where $\Phi_m^*\varphi_Q$ is as in (3.4). 
Then for $m \gg 1$,
\begin{align*}
C_{26} \; &= \; \max_{J \in \Sigma (\frak p )} 
\int_M |{J}|^{\,2}_{\omega_0}\, \omega_0^n\; \geq  \; 
\int_M  |{Q}|^{\,2}_{\omega_0}\; \omega_0^n\; \geq \; 
C_{27} \int_M  |{Q} |^{\,2}_{\eta_m}\; \eta_m^{\,n}\\
&= \; C_{27}\,m\,\int_M \left \{ 
\frac{\Sigma_{\alpha = 0}^{N_m} \,
\hat{\gamma}_{\alpha\alpha}(Q)^2\, |\tilde{s}_{\alpha}|_{h (\ell )}^2}{m^2\Sigma_{\alpha =
0}^{N_m}\,
 |\tilde{s}_{\alpha}|_{h (\ell )}^2}
\,-\, \phi_Q^{\,2} \right \}\,\eta_m^{\,n}\; \\
&\geq\; \frac{C_{28}}{\; m^{n+1}} \,\left \{ \int_M \Sigma_{\alpha = 0}^{N_m} \,
\hat{\gamma}_{\alpha\alpha} (Q)^2\,|\tilde{s}_{\alpha}|_{h (\ell )}^2\,\eta_m^{\,n}\right \} \, -\,
C_{29}\,m
\\ &\geq \; \frac{C_{30}}{\; m^{n+1}} \,\left \{ \int_M \Sigma_{\alpha = 0}^{N_m} \,
\hat{\gamma}_{\alpha\alpha}(Q)^2\, |\tilde{s}_{\alpha}|_{h (\ell )}^2\,\omega (\ell )^{n}\right \} \, -\,
C_{29}\,m
\\ & \geq  \; C_{31} \frac{\| \hat{\gamma}_Q\|_m^{\,2}}{m^{n+1}} \, -\, C_{29}\,m,
\end{align*}
for some $C_{26}$, $C_{27}$, $C_{28}$, $C_{29}$, $C_{30}$ and $C_{31}$.
Hence, if $m \gg 1$, then 
$$
\| \gamma_Q\|_m^{\,2} \; \leq \; 
 \| \hat{\gamma}_Q\|_m^{\,2}\; \leq \; C_{32}\, m^{n+2}.
 \leqno{(3.37)}
$$
for some $C_{32}$. 
Now by $\zeta (0) = 1$, we define a real analytic function $\tilde{\zeta} = \tilde{\zeta}(x) $ on $\Bbb R$ 
satisfying $\tilde{\zeta}(0) = 0$ by
$$
\tilde{\zeta}(x) := \zeta (x)
- 1.
$$
For $X\in \frak p''_m$ above,
by choosing an orthonormal basis for $(V_m, \rho_{h (\ell )})$ possibly distinct from the original one,
we may assume that the representation matrix $\gamma_X $ of $X$ is a real diagonal matrix.
Recall that $\tilde{X}=\, t\,X$, where $|t| < \delta_0 := q^{1/2 + \ell_0 }/\| \hat{\gamma}_{X''} \|_m$.
Put $\tilde{X}'' := t \,{X}''$.
Then by (3.18),
$$
\sqrt{q/C_{11}}\, \| \gamma_{\tilde{X}} \|_m \; \leq \; 
 \| \hat{\gamma}_{\tilde{X}''}\|_m \; =\;  
|t|\cdot \| \hat{\gamma}_{X''} \|_m \; \leq \; q^{1/2 + \ell_0},
$$
i.e., $\| \gamma_{\tilde{X}} \|_m  \, \leq \, \sqrt{C_{11}}\, q^{\ell_0}$. 
Hence, if $m \gg 1$, 
$$
\| \,\gamma^{}_{\tilde{\zeta} (\operatorname{ad} \tilde{X} ) Q}\,\|_m \; \leq \; C_{33} \, q^{\ell_0}
\,\|\gamma^{}_Q\|_m,
 \leqno{(3.38)}
$$
for some $C_{33}$.
Now by the same argument as in (3.36), we see that, for some $C_{34}$,
$$
\|\gamma^{}_{\tilde{\zeta} (\operatorname{ad} \tilde{X} ) Q}\,\|_m^{\,2} \; \geq \; 
C_{34}\,  m^{n+2} \int_M \,\Phi_m^*\varphi^2_{\tilde{\zeta} (\operatorname{ad} \tilde{X} ) Q}
\,\eta_m^{\,n},
\qquad \text{if $m  \gg 1$.}
\leqno{(3.39)}
$$
Put $a_m := \sqrt{\int_M \Phi_m^*\varphi^2_{\tilde{\zeta} (\operatorname{ad} \tilde{X} )
Q}\eta_m^{\,n}}$. Then for $m \gg 1$, by (3.37), (3.38) and (3.39),
$$
a_m \; \leq \; C_{35} \,q^{\ell_0} \leqno{(3.40)}
$$
for some $C_{35}$.
Consider the Laplacians $\square_{\eta_m}$ and $\square_{\omega_0}$
on functions for the K\"ahler manifolds 
$(M, \eta_m )$ and $(M, \omega_0)$, respectively.
Note that $\zeta (\operatorname{ad} \tilde{X} ) Q = Q + \tilde{\zeta} (\operatorname{ad} \tilde{X} ) Q$.
Then for $m \gg 1$,  by (3.35), 
 we obtain
\begin{equation}
\begin{cases}
&
\left | \int_M (\square_{\eta_m}
\phi^{}_Q ) \phi^{}_{\zeta (\operatorname{ad} \tilde{X} ) Q} \eta_m^{\,n}\right |\\
& =\; \left | \int_M   (\square_{\eta_m}
\phi^{}_Q ) \left (    \phi_Q^{}
 + \Phi_m^*\varphi^{}_{\tilde{\zeta} (\operatorname{ad} \tilde{X} ) Q}\right ) 
\,\eta_m^{\,n}\right |  
 \\
& 
\geq \; \| \bar{\partial} \phi_Q \|^{\,2}_{L^2(M, \eta_m )} \; - \; \left | \int_M 
(\square_{\eta_m}
\phi^{}_Q ) 
( \Phi_m^*\varphi^{}_{\tilde{\zeta} (\operatorname{ad} \tilde{X} ) Q})
\,\eta_m^{\,n}\right | \\
&\geq \; \| \bar{\partial}\{ u_Q + q (Q\xi_m )\} \|^{\,2}_{L^2(M, \eta_m )} \\
&\qquad \qquad - a_m \| 
\square_{\eta_m}
\{ u_Q + q (Q\xi_m )\} \|^{}_{L^2 (M, \eta_m )}
\\
&\geq \; (1- \epsilon )\, R_m\;  -\; 
 (1+ \epsilon ) \,a_mS_m,
\end{cases} \tag{3.41}
\end{equation}
where we put $R_m:= \| \bar{\partial}\{ u_Q + q (Q\xi_m )\} \|^{\,2}_{L^2(M, \omega_0 )}$ 
and $S_m := \|  \square_{\omega_0}
\{ u_Q + q (Q\xi_m )\} \|^{}_{L^2 (M, \omega_0 )}$, and
$\epsilon \ll 1$ is a positive constant independent of the choice of $m$, $X$ and $Q$.
Since $Q$ belongs to the compact set $\Sigma (\frak p )$, by (3.34) and the
equality $\|
\bar{\partial}u_Q \|^{}_{L^2(M, \omega_0)} = 1$, we obtain constants $C_{36}$ and $C_{37}$ such that
\begin{equation}
\begin{cases}
&R_m
\; \geq 1\, -\, 2\, q\,\| \bar{\partial}(Q\xi_m ) \|_{L^2(M, \omega_0 )}\; 
\geq \; 1\, - \, C_{36} \,q,\\
&S_m\;
\leq \; \|  \square_{\omega_0} u_Q  \|^{}_{L^2 (M, \omega_0 )}\, +\,
q\,\|  \square_{\omega_0}  (Q\xi_m ) \|^{}_{L^2 (M, \omega_0 )}\; \leq \; C_{37}.
\end{cases}\tag{3.42}
\end{equation}
Then for $m \gg 1$, by (3.40), (3.41) and (3.42), we finally obtain
$$
\left | \int_M (\square_{\eta_m}
\phi^{}_Q )\, \phi^{}_{\zeta (\operatorname{ad} \tilde{X} ) Q} \,\eta_m^{\,n}\right |
\; 
\geq \;   (1- \epsilon )\, (1- C_{36} q) \, 
-\, (1+ \epsilon ) \,C_{35} C_{37}\,q^{\ell_0}\; 
>\; 0,
$$
which implies (3.33), as required.

\end{document}